\documentclass[oneside,english]{amsart}
\usepackage[T1]{fontenc}
\usepackage[latin9]{inputenc}
\usepackage{textcomp}
\usepackage{amstext}
\usepackage{amsthm}
\usepackage{amssymb}

\makeatletter

\providecommand{\tabularnewline}{\\}

\numberwithin{equation}{section}
\numberwithin{figure}{section}
\theoremstyle{plain}
\newtheorem{thm}{\protect\theoremname}
\theoremstyle{definition}
\newtheorem{defn}[thm]{\protect\definitionname}
\theoremstyle{remark}
\newtheorem{rem}[thm]{\protect\remarkname}
\theoremstyle{plain}
\newtheorem{lem}[thm]{\protect\lemmaname}

\makeatother

\usepackage{babel}
\providecommand{\definitionname}{Definition}
\providecommand{\lemmaname}{Lemma}
\providecommand{\remarkname}{Remark}
\providecommand{\theoremname}{Theorem}

\begin{document}
\title{Hodge numbers of arbitrary sections from linear sections}
\author{Herbert Clemens}
\address{Mathematics Dept., Ohio State University, Columbus OH 43210, USA}
\dedicatory{To Enrico Arbarello, my first graduate student and dear life-long
friend.}
\date{November 4, 2022}
\begin{abstract}
Let $\left|L\right|$ be the total space of the inverse of a very
ample line bundle $\pi:L^{-1}\rightarrow B$ over a projective manifold
$B$. Any section of $L^{-1}\rightarrow B$ is isomorphic to $B$
and the Hodge numbers of any proper smooth multisection are determined
by the degree $d$ of that multi-section as are the Hodge numbers
of any smooth complete intersection of multi-sections of degrees $\left(d_{1},\ldots,d_{r}\right)$.
In this paper recursive formulae are given for those Hodge numbers
in terms of the integers $\left\{ d_{1},\ldots,d_{r}\right\} $ and
the Hodge numbers of the linear sections
\[
B,B\text{·}B,\ldots,B^{\dim B}.
\]
The recursion proceeds by induction on dimension and degree. Its proof
relies on the theory of asymptotic mixed Hodge structures. 

An interesting corollary is that the Lefschetz hyperplane property
is weakened by one degree in this setting. That is, relative vanishing
does not reach the middle degree of the hyperplane section but only
to degree one less than the middle degree.

As an application, in an Appendix, we calculate closed formulae for
all Hodge numbers of all smooth complete intersections for the case
$\dim B=3$.
\end{abstract}

\email{clemens.43@osu.edu}

\maketitle
\tableofcontents{}

\section{Introduction}

Long ago the author was among those who worked out the asymptotic
mixed Hodge theory for a one-parameter degeneration of Kähler manifolds
in the complex analytic setting \cite{Clemens-0}. Quite recently,
as part of a collaboration with string theorists, I was called upon
to compute the Hodge numbers of certain elliptically fibered Calabi-Yau
threefolds and fourfolds in settings where standard computational
tools are less accessible. Even though it is well established that
such formulae 'can be computed,' I found it annoying that I could
not find in the literature closed formulae for Hodge numbers of complete
intersections of a fixed complex projective manifold of dimension
$n+1$ whose Hodge numbers were given, or more generally of complete
intersections of multi-sections of an ample line bundle over a complex
projective manifold of dimension $n$whose Hodge numbers were given,
even for low values of $n$. The purpose of this paper is lay out
the recursive algorithm that allows the derivation of such formulas
in all dimensions and to derive such formulas for $n\leq3$. 

Although the computations involved may require some patience to follow,
the principle by which they are derived is quite simple, any complete
intersection $V$ specializes linearly to the transverse union of
two complete intersections $V'\cup V''$ of lower degree unless the
degree was one to start with. Then asymptotic mixed Hodge theory allows
one to read off the Hodge numbers of the original complete intersection
from those of the two components $V'$ and $V''$ of the specialization,
the intersection $V'\cap V''$ of those two components, and finally
the intersection $V\cap V'\cap V''$ of $V$ with the intersection
of those two components. 

The computational method is by complete induction on degrees and dimension,
using asymptotic mixed Hodge theory 'in reverse.' Using modern mathematical
software it is presumably a simple matter to write a closed formula
for the Hodge numbers of all complete intersections for any fixed
$n$. In the Appendix to this paper we provide a roadmap for such
a program by working out the formulas for all proper complete intersections
in the case in which $n\leq3.$

'In reverse' simply means that since the dimensions of the graded
pieces of the Hodge filtration of the asymptotic mixed Hodge structure
are the same as those for the Hodge structure of the nearby smooth
fiber (see \cite{Peters}, §13), and the latter can be more easily
computed by computing the former, namely for a degeneration that breaks
the nearby fiber into simpler pieces in the limit as desribed above.
This strategy is certainly not a new one, going way back to the 19th-century
technique of computing the genus of a smooth plane curve of degree
$d$ by degenerating it into a union of $d$ projective lines meeting
transversally, or inductively, degenerating it into a union of two
transversely intersecting smooth curves of degrees $d_{1}$ and $d_{2}$
respectively where $d_{1}+d_{2}=d$. Our setting in Case One below
is projective submanifolds of the total space of $Y=\left|L\right|$
of the inverse of a very ample line bundle $\pi:L^{-1}\rightarrow B$
over a projective manifold $B$. The zero scheme of any section of
$L^{-1}\rightarrow B$ is isomorphic to $B$ and the Hodge numbers
of any proper smooth multisection are determined by the degree $d$
of that multi-section as are the Hodge numbers of any smooth complete
intersection of multi-sections of degrees $\left(d_{1},\ldots,d_{r}\right)$.
In this paper recursive formulae are given for computing those Hodge
numbers in terms of the integers $\left\{ d_{1},\ldots,d_{r}\right\} $
and the Hodge numbers of the linear sections. The same method works
when $Y$ is a project manifold (Case Two below) although in that
case, the derivation of Hodge numbers of complete intersections is
well known.

However we believe the case where $Y$ is the total space of a negative
line bundle on a projective manifold (Case One below) is new. We include
both cases in the combined exposition below, since it requires no
additional work and the explicit formulae in the case in which $Y$
is compact are needed in the derivation in which $Y$ is non-compact.
Also, it seems that these formulas, even in the $Y$-compact case,
do not explicitly appear in the literature. In short each induction
step is always established as described above, namely by using that,
if a smooth proper divisor $V_{d}$ specializes linearly to a sum
of smooth divisors $V_{d_{1}}+V_{d_{2}}$ where $d=d_{1}+d_{2}$,
allowing us to derive the Hodge numbers $h^{p,q}\left(V_{d}\right)$
as simple sums of certain Hodge numbers of $V_{d_{1}},$ $V_{d_{2}}$,
$V_{d_{1}}\cap V_{d_{2}}$, and $V_{d}\cap V_{d_{1}}\cap V_{d_{2}}$. 

A corollary in Case One of the derivation is that the Lefschetz hyperplane
property is weakened by one degree in the Case One setting. That is,
relative vanishing does not reach the middle degree of the hypersurface
section but only to degree one less than the middle degree.

\section{Computing Hodge numbers by induction on dimension and degree}

\subsection{Case One: Projective submanifolds of the total space of a negative
line bundle on a projective manifold}

Let
\[
B\subseteq\mathbb{P}^{M}=\mathbb{P}\left(H^{0}\left(L^{-1}\right)\right)
\]
be a smooth complex projective manifold of dimension $n$ imbedded
by the complete linear system of a very ample line bundle $L/B$.
We assume that $B$ does not contain the point $\left[1,0,\ldots,0\right]\in\mathbb{P}^{M}$.Then
the zero-schemes of homogeneous forms 
\[
F\left(y\right)\in\mathrm{Sym}{}^{d}\left(H^{0}\left(L\right)\right)
\]
be a line bundle such that $L^{-1}=\mathcal{O}_{B}\left(N\right)$
is very ample so that its complete linear system gives an imbedding
\[
B\rightarrow\mathbb{P}\left(H^{0}\left(L^{-1}\right)\right)=\mathbb{P}^{M},
\]

We denote homogeneous coordinates
\[
\mathbb{P}_{\left[y_{0},\ldots,y_{M}\right]}=\mathbb{P}^{M}
\]
and let $\left|\mathcal{O}_{\mathbb{P}^{M}}\left(-1\right)\right|$
denote the total space of the line bundle $\mathcal{O}_{\mathbb{P}^{M}}\left(-1\right)$
and form the fibered product
\[
\begin{array}{ccc}
Y & \rightarrow & \left|\mathcal{O}_{\mathbb{P}^{M}}\left(-1\right)\right|\\
\downarrow &  & \downarrow\\
B & \rightarrow & \mathbb{P}^{M}.
\end{array}
\]

\begin{defn}
We call a homogeneous form  
\[
F\left(y\right)\in\mathrm{Sym}{}^{d}\left(H^{0}\left(L\right)\right)=\mathrm{Sym}{}^{d}\left(\mathbb{C}\left[y_{0},\ldots,y_{M}\right]\right)
\]
\textit{monic} if $F\left(\left[1,0,\ldots,0\right]\right)\neq0$,
that is, the coefficient of $y_{0}^{d}$ in $F\in\mathbb{C}\left[y_{0},\ldots,y_{M}\right]$
is not zero.
\end{defn}

If $F\left(y\right)$ is monic, then we have 'affine' coordinate $\vartheta=y_{0}\in H^{0}\left(\mathcal{O}_{\mathbb{P}\left(H^{0}\left(L^{-1}\right)\right)}\left(1\right)\right)$
and we can write
\[
F\left(y\right)=\sum_{j=0}^{d}f_{j}\left(y_{1},\ldots,y_{M}\right)\cdot\vartheta^{d-j}where
\]
$f_{j}$ is homogeneous of degree $j$. Also 
\[
\left(F\right):=\left\{ y:F\left(y\right)=0\right\} \subseteq Y
\]
is proper over $B$ and forms a $d$-sheeted branched cover of $B$.
Given monic $F_{d}\left(y\right)$ and $G_{d}\left(y\right)$ we can
form a linear family
\[
\left\{ t\text{·}F_{d}\left(y\right)+G_{d}\left(y\right)=0\right\} \subseteq\Delta\times Y,
\]
so that all of its forms for $\left|t\right|<\varepsilon$ are monic
and so its fibers over the $t$-disk $\Delta$ are proper. Thus 'linear
variation' and' degeneration of projective varieties' makes sense
in this context. For monic $F_{1}\left(y\right)\in H^{0}\left(\mathcal{O}_{\mathbb{P}^{M}}\left(1\right)\right)$
and $\left\{ F_{1}\left(y\right)=0\right\} \cap Y=:V_{1}$ is isomorphic
to $B$, $V_{1}\text{·}V_{1}$ is a divisor in the linear system associated
with $L^{-1}$, etc., and finally the cardinality of $V_{1}^{\dim B}$
is the degree of $B\subseteq\mathbb{P}^{M}$.

For each degree $d$ and monic $F_{d}$ we write
\[
V_{d}:=\left\{ F_{d}=0\right\} \subseteq\left|L\right|
\]
and if $d=d_{1}+d_{2}$, we can form the linear degeneration
\begin{equation}
\left\{ t\text{·}F_{d}\left(y\right)-F_{d_{1}}\left(y\right)\text{·}F_{d_{2}}\left(y\right)=0\right\} \subseteq\mathbb{C}\times\left|L\right|.\label{eq:family0}
\end{equation}
Our goal will be to explicitly derive the Hodge numbers of 
\[
V_{d_{1}}\cdot\ldots\cdot V_{d_{r}}
\]
for $r<n$ from the Hodge numbers of $B$. Of course if $r=n$, the
zero-th Hodge number is 
\[
d_{1}\cdot\ldots\cdot d_{n}\cdot\deg B.
\]

\subsection{Case Two: Smooth hypersurface sections of a projective manifold}

An alternative situation is that in which $Y$ is a projective manifold
of dimension $n+1$ and $B$ is a smooth hyperplane section. Here
the 'monic' condition is irrelevant but the same recursive formula
given in Theorem \ref{thm:In-either-Case} below applies. In this
case one uses the same reasoning as in Case One to derive explicit
recursive formulas for the Hodge numbers of complete intersections
of any projective manifold from the Hodge numbers of the manifold
$Y$ and the (middle) Hodge numbers of a smooth hyperplane section.
In this case Theorem \ref{thm:In-either-Case} below applies to the
setting 
\begin{equation}
\left\{ t\text{·}F_{d}\left(y\right)-F_{d_{1}}\left(y\right)\text{·}F_{d_{2}}\left(y\right)=0\right\} \subseteq\mathbb{C}\times Y\subseteq\mathbb{C}\times\mathbb{P}^{M}\label{eq:family1}
\end{equation}
where $B$ is a smooth hyperplane section of $Y$, $\mathcal{O}_{B}\left(N\right)=\left.\mathcal{O}_{\mathbb{P}^{M}}\left(1\right)\right|_{B}$,
and $F_{d}$ is any homogeneous form of degree defining a smooth hypersurface
of $Y$. 

\subsection{Recursive formula}

In either Case One or Case Two, let 
\[
V_{d_{1}}\text{·}\ldots\text{·}V_{d_{r}}
\]
denote the transverse intersection of $r$ smooth proper hypersurfaces
of $Y$of degrees $d_{1},\ldots,d_{r}$ respectively. The purpose
of note is to prove the following:
\begin{thm}
\label{thm:In-either-Case}In either Case One or Case Two above, given
the Hodge numbers of linear sections
\[
V_{1},V_{1}^{2},\ldots,V_{1}^{n},
\]
the Hodge numbers of $V_{d_{1}}\text{·}\ldots\text{·}V_{d_{r}}$ can
be computed recursively by complete induction on dimension and degree
using only the following formulae:

Let $V_{d}\Rightarrow V_{d_{1}}+V_{d_{2}}$ be as in either (\ref{eq:family0})
or (\ref{eq:family1}) above , and let $\tilde{V}_{d_{2}}$ denote
the blow-up of $V_{d_{2}}$ blown up along the submanifold $V_{d}\cap V_{d_{1}}\cap V_{d_{2}}$. 

\noindent For $p+q=k>n+1$
\[
h^{p,q}\left(V_{d}\right)=h^{p,q}\left(V_{1}\right),
\]
for $k=p+q=n+1$
\[
\begin{array}{c}
h^{p,q}\left(V_{d}\right)=h^{p,q}\left(\ker\left(H^{n+1}\left(V_{d_{1}}\right)\oplus H^{n+1}\left(\tilde{V}_{d_{2}}\right)\rightarrow H^{n+1}\left(V_{d_{1}}\cap V_{d_{2}}\right)\right)\right)-h^{p,q}\left(V_{d_{1}}\cap V_{d_{2}}\right)\\
=\left\{ \begin{array}{c}
\dim\left(\ker\left(H^{p,q}\left(V_{d_{1}}\right)\oplus H^{p,q}\left(V_{d_{2}}\right)\rightarrow H^{p,q}\left(V_{d_{1}}\cap V_{d_{2}}\right)\right)\right)\\
+h^{p-1,q-1}\left(V_{d}\cap V_{d_{1}}\cap V_{d_{2}}\right)-h^{p,q}\left(V_{d_{1}}\cap V_{d_{2}}\right)
\end{array}\right.,
\end{array}
\]
and for $p+q=k=n$,
\[
h^{p,q}\left(V_{d}\right)=\left\{ \begin{array}{c}
h_{prim}^{p-1,q}\left(V_{d_{1}}\cap V_{d_{2}}\right)\\
h_{prim}^{p,q}\left(V_{d_{1}}\right)+h_{prim}^{p,q}\left(V_{d_{2}}\right)+h^{p-1,q-1}\left(V_{d}\cap V_{d_{1}}\cap V_{d_{2}}\right)\\
+h_{prim}^{p,q-1}\left(V_{d_{1}}\cap V_{d_{2}}\right).
\end{array}\right.
\]
\end{thm}

\begin{rem}
The complete induction in the above theorem in Case One requires application
of either Case One or Case Two with prior intersections playing the
role of $Y$ at various steps.
\end{rem}

As an example of the application of Case Two of the formulae in Theorem
\ref{thm:In-either-Case}, we calculate the Hodge numbers of the quintic
threefold $X_{5}\subseteq\mathbb{P}^{4}=Y$ by degenerating it linearly
\[
\left\{ t\text{·}G_{5}+G_{3}\text{·}G_{2}=0\right\} \subseteq\mathbb{C}\times\mathbb{P}^{4}
\]
into the union of a cubic threefold $X_{3}$ and a quadric threefold
$X_{2}$, taking care that the intersection $X_{5}\cap X_{3}\cap X_{2}$
is transverse. Rewriting the equation of the family as
\[
\left\{ \left|\begin{array}{cc}
t & G_{2}\\
G_{3} & G_{5}
\end{array}\right|=0\right\} \subseteq\mathbb{C}\times\mathbb{P}^{4}
\]
we must do a small resolution over the curve $C:=\left\{ t=0\right\} \cap X_{5}\cap X_{3}\cap X_{2}$
where all four entries of the above matrix vanish. This has the effect
of blowing up $C$ in either $X_{2}$ or $X_{3}$ and while modlifying
the degeneration to normal crossing form without base locus as considered
in \cite{Clemens-0}. Noting that $X_{2}\cap X_{3}$ is a $K3$-surface,
Theorem \ref{thm:In-either-Case} and the table\smallskip{}
\[
\begin{array}{ccccc}
 & X_{3} & X_{2} & X_{3}\cap X_{2} & X_{5}\cap X_{3}\cap X_{2}\\
h^{3,0} & 0 & 0\\
h^{2,1} & 5 & 0\\
h^{2,0} & 0 & 0 & 1\\
h_{prim}^{1,1} & 0 & 0 & 19\\
h^{1,0} &  &  &  & 76
\end{array}
\]
allow us to compute conclude that
\[
\begin{array}{c}
h^{3,0}\left(X_{5}\right)=h^{3,0}\left(X_{3}\right)+h^{3,0}\left(X_{2}\right)+h^{2,0}\left(X_{3}\cap X_{2}\right)=1\\
h^{2,1}\left(X_{5}\right)=\left(\begin{array}{c}
h^{2,1}\left(X_{3}\right)+h^{2,1}\left(X_{2}\right)+h^{2,0}\left(X_{3}\cap X_{2}\right)+\\
h_{prim}^{1,1}\left(X_{3}\cap X_{2}\right)+h^{1,0}\left(X_{5}\cap X_{3}\cap X_{2}\right)
\end{array}\right)=101.
\end{array}
\]

\section{Common setting for Cases One and Two}

\subsection{Small resolution of linear family over the $t$-disk}

We begin by rewriting (\ref{eq:family0}) or (\ref{eq:family1}) in
the form
\[
\left\{ \left|\begin{array}{cc}
t & F_{d_{1}}\\
F_{d_{2}} & F_{d}
\end{array}\right|=0\right\} 
\]
from which one sees that the total space of the family in $\mathbb{C}\times\mathbb{P}^{M}$
has a nodal locus on the set
\[
\left\{ t=F_{d_{1}}=F_{d_{2}}=F_{d}=0\right\} .
\]

Using either ratios of rows or ratios of columns in the above matrix
gives two small resolutions of the total space of the family and deposits
the exceptional locus in either $V_{d_{1}}$ or $V_{d_{2}}$. We choose
$V_{d_{2}}$ which is thereby blown up along the codimension-two submanifold
$\left\{ F_{d_{1}}=F_{d}=0\right\} .$ We denote the blown up $V_{d_{2}}$
as $\tilde{V}_{d_{2}}$. We have Hodge subspaces
\[
H^{j,k}\left(\tilde{V}_{d_{2}}\right)=H^{j,k}\left(V_{d_{2}}\right)\oplus H^{j-1,k-1}\left(V_{d_{1}}\cap V_{d_{2}}\cap V_{d}\right)\otimes\left(1,1\right)
\]
where $\left(1,1\right)$ denotes the so-called Hodge-Tate mixed Hodge
structure (dimension $1$, weight $2$, pure type $\left(1,1\right)$).
Furthermore\cite{Clemens-0}
\[
V_{d_{1}}\cap\tilde{V}_{d_{2}}=V_{d_{1}}\cap V_{d_{2}}.
\]
A detailed description of this normal-crossing degeneration and its
asymptotic Hodge theory is given in \cite{Clemens-0}. 

\subsection{Topological decomposition of $V_{d}$}

$V_{d}$ is constructed topologically by first removing a small regular
open neighborhood $U_{V_{d_{1}}\cap V_{d_{2}}}$ of $V_{d_{1}}\cap V_{d_{2}}$
from $V_{d_{1}}\cup\tilde{V}_{d_{2}}$, giving respectively $V'_{d_{1}}$
as a deformation retraction of $V_{d_{1}}-\left(V_{d_{1}}\cap V_{d_{2}}\right)$
and $\tilde{V}'_{d_{2}}$ as a deformation retraction of $\tilde{V}_{d_{2}}-\left(V_{d_{1}}\cap V_{d_{2}}\right)$.
Then $\partial V'_{d_{1}}\cong\partial\tilde{V}'_{d_{2}}=:T$ where
$T$ is a circle bundle over $\left(V_{d_{1}}\cap V_{d_{2}}\right)$
and the isomorphism reverses orientation. $V_{d}$ is then given topologically
by pasting $V'_{d_{1}}$ to $\tilde{V}'_{d_{2}}$ by identifying corresponding
points of their common boundary. There is a natural 'contraction'
mapping
\begin{equation}
\rho:V_{d}\rightarrow V_{d_{1}}\cup\tilde{V}_{d_{2}}\label{eq:contract}
\end{equation}
obtained by the inclusion of $V_{d}$ into a regular neighborhood
of $V_{d_{1}}\cup\tilde{V}_{d_{2}}$ in the (smooth) total space of
the family. These constructions are described in detail in §5-7 of
\cite{Clemens-0} where the topology of the horizontally and vertically
exact diagram
\noindent \begin{flushleft}
\begin{equation}
\begin{array}{ccccc}
\ldots &  & \ldots &  & \ldots\\
\downarrow &  & \downarrow &  & \downarrow\\
\ldots\rightarrow H^{\cdot}\left(V_{d_{1}}\cup\tilde{V}_{d_{2}}\right) & \rightarrow & H^{\cdot}\left(V_{d_{1}}\right)\oplus H^{\cdot}\left(\tilde{V}_{d_{2}}\right) & \rightarrow & H^{\cdot}\left(V_{d_{1}}\cap V_{d_{2}}\right)\rightarrow\ldots\\
\downarrow^{\rho^{\ast}} &  & \downarrow &  & \downarrow^{\rho^{\ast}}\\
\ldots\rightarrow H^{\cdot}\left(V_{d}\right) & \rightarrow & \begin{array}{c}
H^{\text{·}}\left(A_{V_{d_{1}}}^{\text{·}}\left(\log\left(V_{d_{1}}\cap V_{d_{2}}\right)\right)\right)\oplus\\
H^{\text{·}}\left(A_{\tilde{V}{}_{d_{2}}}^{\text{·}}\left(\log\left(V_{d_{1}}\cap V_{d_{2}}\right)\right)\right)
\end{array} & \rightarrow & H^{\text{·}}\left(T\right)\rightarrow\ldots\\
\downarrow &  & \downarrow^{\mathrm{Residue}} &  & \downarrow^{\mathrm{Residue}}\\
\ldots\rightarrow\mathbb{H}^{\text{·}+1}\left(\rho^{\ast}\right) & \rightarrow & \begin{array}{c}
H^{\text{·}-1}\left(V_{d_{1}}\cap V_{d_{2}}\right)\otimes\left(1,1\right)\oplus\\
H^{\text{·}-1}\left(V_{d_{1}}\cap V_{d_{2}}\right)\otimes\left(1,1\right)
\end{array} & \rightarrow & H^{\cdot-1}\left(V_{d_{1}}\cap V_{d_{2}}\right)\otimes\left(1,1\right)\rightarrow\ldots\\
\downarrow &  & \downarrow^{\mathrm{push-forward}} &  & \downarrow^{\text{·}d_{1}N}\\
\ldots\rightarrow H^{\cdot+1}\left(V_{d_{1}}\cup\tilde{V}_{d_{2}}\right) & \rightarrow & H^{\cdot+1}\left(V_{d_{1}}\right)\oplus H^{\cdot+1}\left(\tilde{V}_{d_{2}}\right) & \rightarrow & H^{\cdot+1}\left(V_{d_{1}}\cap V_{d_{2}}\right)\rightarrow\ldots\\
\downarrow^{\rho^{\ast}} &  & \downarrow &  & \downarrow^{\rho^{\ast}}\\
\ldots\rightarrow H^{\cdot+1}\left(V_{d}\right) & \rightarrow & \begin{array}{c}
H^{\text{·}+1}\left(A_{V_{d_{1}}}^{\text{·}}\left(\log\left(V_{d_{1}}\cap V_{d_{2}}\right)\right)\right)\oplus\\
H^{\text{·}+1}\left(A_{\tilde{V}{}_{d_{2}}}^{\text{·}}\left(\log\left(V_{d_{1}}\cap V_{d_{2}}\right)\right)\right)
\end{array} & \rightarrow & H^{\text{·}+1}\left(T\right)\rightarrow\ldots\\
\downarrow &  & \downarrow^{\mathrm{Residue}} &  & \downarrow^{\mathrm{Residue}}\\
\ldots\rightarrow\mathbb{H}^{\text{·}+2}\left(\rho^{\ast}\right) & \rightarrow & \begin{array}{c}
H^{\text{·}}\left(V_{d_{1}}\cap V_{d_{2}}\right)\otimes\left(1,1\right)\oplus\\
H^{\text{·}}\left(V_{d_{1}}\cap V_{d_{2}}\right)\otimes\left(1,1\right)
\end{array} & \rightarrow & H^{\cdot}\left(V_{d_{1}}\cap V_{d_{2}}\right)\otimes\left(1,1\right)\rightarrow\ldots\\
\downarrow &  & \downarrow &  & \downarrow^{\text{·}d_{1}N}\\
\ldots &  & \ldots &  & \ldots
\end{array}\label{eq:biggest}
\end{equation}
is also explained as are the hypercohomology groups $\mathbb{H}^{\text{·}}\left(\rho^{\ast}\right)$
characterized by the fact that their addition completes the cohomology
group mappings
\[
\rho^{\ast}:V_{d}H^{\cdot}\left(V_{d_{1}}\cup\tilde{V}_{d_{2}}\right)\rightarrow H^{\cdot}\left(V_{d}\right)
\]
to a long exact cohomology sequence. The exact cohomology diagram
(\ref{eq:biggest}) can be thought of as intertwining horizontal Mayer-Vietoris
exact sequences with a vertical mapping cone exact sequence and vertical
`log' exact sequences derived from residue maps via the Gysin isomorphism.
\par\end{flushleft}

\subsection{The governing exact diagram of morphisms of mixed Hodge structures}

We now apply the theory of mixed Hodge structures to the topological
set-up described in the previous Subsection, in particular to the
diagram (\ref{eq:biggest}). For readers unfamiliar with mixed Hodge
theory, a gentle introduction can be found in \cite{Durfee}, a more
in-depth treatment in \cite{Griffiths} and finally the full theory
in \cite{Peters}

The cohomology mappings (\ref{eq:biggest}) associated to the mapping
cone of $\rho$ in (\ref{eq:contract}) in the case of the the topological
decomposition of $V_{d}$ give rise to morphisms of mixed Hodge structures
explained in Chapter 3 of \cite{Clemens-0}. Our strategy will be
to deduce the Hodge numbers
\[
H^{\cdot}\left(V_{d};\mathbb{C}\right)=\sum_{p+q=\cdot}H^{p,q}\left(V_{d}\right)
\]
from the Hodge numbers of all the other cohomology groups in (\ref{eq:biggest}),
all of which involve only $H^{\cdot}\left(V'_{d'}\right)$ with $d'<d$
and/or $\dim V'_{d'}<\dim V{}_{d}$. This and closed expressions for
$\sum_{k=1}^{n}k^{m}$ with $m\leq n+1$ then allow us to recursively
build closed formulas for Hodge numbers of proper complete intersections
of $Y$ for $\dim Y\leq n+1$. We carry out the recursive algorithm
in detail for $n\leq3$ in the Appendix to this paper. It seems to
the author that, following the paradigm in the Appendix, computer
programs that make the corresponding calculations for arbitrary (fixed)
$n$ ought not to be too difficult to design. 

\subsection{\label{subsec:Mixed-Hodge-structures}Mixed Hodge structures}

For all the mixed Hodge structures we will consider, $W_{\cdot}$
will denote the weight filtration and $F^{\cdot}$ will denote the
Hodge filtration. Taking to account that the rank-one Hodge-Tate structure
$\left(1,1\right)$ has weight $2$, all weight and Hodge filtration
of terms are just the the standard ones of the cohomology groups represented,
except for $H^{\text{·}}\left(T\right)$. The right-hand vertical
exact sequence

\[
\ldots\overset{\text{·}d_{1}N}{\longrightarrow}H^{k}\left(V_{d_{1}}\cap V_{d_{2}}\right)\overset{\rho^{\ast}}{\longrightarrow}H^{k}\left(T\right)\overset{\mathrm{Residue}}{\longrightarrow}H^{k-1}\left(V_{d_{1}}\cap V_{d_{2}}\right)\overset{\text{·}d_{1}N}{\longrightarrow}H^{k+1}\left(V_{d_{1}}\cap V_{d_{2}}\right)\overset{\rho^{\ast}}{\longrightarrow}\ldots
\]
is the standard exact sequence for the circle bundle $T/\left(V_{d_{1}}\cap V_{d_{2}}\right)$.
$A^{\text{·}}\left(T\right)$ is quasi-isomorphic to a mapping cone
of
\[
\left(d_{1}N\text{·}\right):H^{\cdot-2}\left(V_{d_{1}}\cap V_{d_{2}}\right)\otimes\left(1,1\right)\rightarrow H^{\cdot}\left(V_{d_{1}}\cap V_{d_{2}}\right)
\]
in the abelian category of mixed Hodge structures and that therefore
endows $H^{k}\left(T\right)$ with its mixed Hodge structure 
\[
\begin{array}{c}
\begin{array}{c}
W_{k}\left(H^{k}\left(T\right)\right)=\mathrm{image}\left(\rho^{\ast}\right)\\
=\mathrm{coker}\left(\left(d_{1}N\text{·}\right):H^{k-2}\left(V_{d_{1}}\cap V_{d_{2}}\right)\otimes\left(1,1\right)\rightarrow H^{k}\left(V_{d_{1}}\cap V_{d_{2}}\right)\right)
\end{array}\\
\frac{W_{k+1}\left(H^{k}\left(T\right)\right)}{W_{k}\left(H^{k}\left(T\right)\right)}=\ker\left(\left(d_{1}N\text{·}\right):H^{k-1}\left(V_{d_{1}}\cap V_{d_{2}}\right)\otimes\left(1,1\right)\rightarrow H^{k+1}\left(V_{d_{1}}\cap V_{d_{2}}\right)\right).
\end{array}
\]
Since $N$ is ample
\[
\ker\left(\left(d_{1}N\text{·}\right):H^{k-1}\left(V_{d_{1}}\cap V_{d_{2}}\right)\otimes\left(1,1\right)\rightarrow H^{k+1}\left(V_{d_{1}}\cap V_{d_{2}}\right)\right)=0
\]
if $k<n$ and 
\[
\mathrm{coker}\left(\left(d_{1}N\text{·}\right):H^{k-2}\left(V_{d_{1}}\cap V_{d_{2}}\right)\otimes\left(1,1\right)\rightarrow H^{k}\left(V_{d_{1}}\cap V_{d_{2}}\right)\right)=0
\]
if $k\geq n$. So
\begin{equation}
\begin{array}{c}
H^{k}\left(T\right)=W_{k}\left(H^{k}\left(T\right)\right)\,\,\,if\,k<n\\
H^{k}\left(T\right)=\frac{W_{k+1}\left(H^{k}\left(T\right)\right)}{W_{k}\left(H^{k}\left(T\right)\right)}\,\,\,if\,k\geq n.
\end{array}\label{eq:weights2}
\end{equation}
Furthermore
\begin{equation}
\frac{H^{k}\left(V_{d_{1}}\cap V_{d_{2}}\right)}{H_{prim}^{k}\left(V_{d_{1}}\cap V_{d_{2}}\right)}\cong\left\{ \begin{array}{c}
H^{k-2}\left(V_{d_{1}}\cap V_{d_{2}}\right)\otimes\left(1,1\right)\,\,if\,k\leq n-1\\
H^{k+1}\left(V_{d_{1}}\cap V_{d_{2}}\right)\,\,if\,k\geq n.
\end{array}\right.\label{eq:weights}
\end{equation}

Also $A^{\text{·}}\left(T\right)$ is quasi-isomorphic to a mapping
cone of
\[
A^{\cdot+1}\left(V_{d}\right)\rightarrow A_{V_{d_{1}}}^{\text{·}+1}\left(\log\left(V_{d_{1}}\cap V_{d_{2}}\right)\right)\oplus A_{\tilde{V}{}_{d_{2}}}^{\text{·}+1}\left(\log\left(V_{d_{1}}\cap V_{d_{2}}\right)\right)
\]
inducing the same mixed Hodge structure.

The cohomology of $V_{d_{1}}\cup\tilde{V}_{d_{2}}$ and its mixed
Hodge structure are given by the isomorphism
\[
H^{k}\left(V_{d_{1}}\cup\tilde{V}_{d_{2}}\right)\cong\mathbb{H}^{k}\left(A^{\cdot}\left(V_{d_{1}}\right)\oplus A^{\cdot}\left(\tilde{V}_{d_{2}}\right)\rightarrow A^{\cdot}\left(V_{d_{1}}\cap\tilde{V}_{d_{2}}\right)\right)
\]
with weight filtration
\[
\begin{array}{c}
W_{k-1}\left(H^{k}\left(V_{d_{1}}\cup\tilde{V}_{d_{2}}\right)\right)=\mathrm{coker}\left(H^{k-1}\left(V_{d_{1}}\right)\oplus H^{k-1}\left(\tilde{V}_{d_{2}}\right)\rightarrow H^{k-1}\left(V_{d_{1}}\cap\tilde{V}_{d_{2}}\right)\right)\\
\frac{H^{k}\left(V_{d_{1}}\cup\tilde{V}_{d_{2}}\right)}{W_{k-1}\left(H^{k}\left(V_{d_{1}}\cup\tilde{V}_{d_{2}}\right)\right)}=\ker\left(H^{k}\left(V_{d_{1}}\right)\oplus H^{k}\left(\tilde{V}_{d_{2}}\right)\rightarrow H^{k}\left(V_{d_{1}}\cup\tilde{V}_{d_{2}}\right)\right).
\end{array}
\]

The mixed Hodge structure on $H^{\cdot}\left(V_{d}\right)$, called
the asymptotic mixed Hodge structure, is explained in \cite{Clemens-0}.
It is given by the isomorphism
\[
H^{\cdot}\left(V_{d}\right)\cong\mathbb{H}^{\cdot}\left(A_{V_{d_{1}}}^{\text{·}}\left(\log\left(V_{d_{1}}\cap V_{d_{2}}\right)\right)\oplus A_{\tilde{V}{}_{d_{2}}}^{\text{·}}\left(\log\left(V_{d_{1}}\cap V_{d_{2}}\right)\rightarrow A^{\text{·}}\left(T\right)\right)\right).
\]
Its weight filtration is given by
\[
\begin{array}{c}
W_{k-1}\left(H^{k}\left(V_{d}\right)\right)=\mathrm{image}\left(W_{k-1}\left(H^{k-1}\left(T\right)\right)\right)\\
W_{k}\left(H^{k}\left(V_{d}\right)\right)=\mathrm{image}\left(H^{k}\left(V_{d_{1}}\cup\tilde{V}_{d_{2}}\right)\right)\\
W_{k+1}\left(H^{k}\left(V_{d}\right)\right)=H^{k}\left(V_{d}\right).
\end{array}
\]
The mapping $\rho^{\ast}$ in (\ref{eq:biggest}) is the mapping $\mu$
in diagram (3.6) of \cite{Clemens-0}. 

The middle vertical sequence is the standard `log' exact sequence
of mixed Hodge structures, and the left-side horizontal maps are morphisms
of mixed Hodge structures. 

\section{The computation}

\subsection{Duality and the Lefschetz hyperplane theorem}

The exact duality of standard exact sequences permit the definitions
\[
\begin{array}{ccc}
\downarrow &  & \uparrow\\
H^{k}\left(V_{d_{1}}\right) &  & H^{2n-k}\left(V_{d_{1}}\right)\\
\downarrow &  & \uparrow^{\mathrm{Push-forward}}\\
H^{k}\left(V_{d_{1}}\cap V_{d_{2}}\right) &  & H^{2n-k-2}\left(V_{d_{1}}\cap V_{d_{2}}\right)\\
\downarrow &  & \uparrow^{\mathrm{Residue}}\\
H^{k+1}\left(V_{d_{1}},\left(V_{d_{1}}\cap V_{d_{2}}\right)\right) &  & H^{2n-\left(k+1\right)}\left(V_{d_{1}}-\left(V_{d_{1}}\cap V_{d_{2}}\right)\right)\\
\downarrow &  & \uparrow\\
H^{k+1}\left(V_{d_{1}}\right) &  & H^{2n-\left(k+1\right)}\left(V_{d_{1}}\right)\\
\downarrow &  & \uparrow^{\mathrm{Push-forward}}
\end{array}
\]
Since $H^{k}\left(V_{d_{1}},\left(V_{d_{1}}\cap V_{d_{2}}\right)\right)=0$
for $k<n-1$ and by duality and ampleness
\[
H^{k}\left(V_{d_{1}}-\left(V_{d_{1}}\cap V_{d_{2}}\right)\right)=0
\]
for $k>n$. Also for $k>n$ we have
\[
\left(H^{2n-2-k}\left(V_{d_{1}}\cap V_{d_{2}}\right)\right)^{\vee}\cong H^{k}\left(V_{d_{1}}\cap V_{d_{2}}\right)\cong H^{k}\left(V_{d_{1}}\right)\cong\left(H^{2n-k}\left(V_{d_{1}}\right)\right)^{\vee}
\]
so that the push-forward map
\[
H^{k}\left(V_{d_{1}}\cap V_{d_{2}}\right)\rightarrow H^{k+2}\left(V_{d_{1}}\right)
\]
is an isomorphism for $k>n-1$ and surjective for $k=n-1$. Then the
commutative diagram
\begin{equation}
\begin{array}{ccc}
H^{k}\left(V_{d_{1}}\cap V_{d_{2}}\right) & \rightarrow & H^{k+2}\left(V_{d_{1}}\right)\\
\downarrow^{\text{·}N} &  & \downarrow^{\text{·}N}\\
H^{k+2}\left(V_{d_{1}}\cap V_{d_{2}}\right) & \rightarrow & H^{k+4}\left(V_{d_{1}}\right)
\end{array}\label{eq:comdia}
\end{equation}
implies that
\begin{equation}
H_{prim}^{k}\left(V_{d_{1}}\cap V_{d_{2}}\right)\rightarrow H_{prim}^{k+2}\left(V_{d_{1}}\right)\label{eq:push}
\end{equation}
induced by push-forward is an isomorphism for $k>n-1$ and surjective
for $k=n-1$. In fact for $k=n-2$ the cokernel of (\ref{eq:push})
is $H_{prim}^{n}\left(V_{d}\right)$ because of the isomorphisms indicated
in the diagram
\begin{equation}
\begin{array}{ccc}
H^{n-4}\left(V_{d_{1}}\cap V_{d_{2}}\right) & \rightarrow & H^{n-2}\left(V_{d_{1}}\right)\\
\downarrow^{\text{·}N} &  & \downarrow^{\text{·}N}\\
H^{n-2}\left(V_{d_{1}}\cap V_{d_{2}}\right) & \rightarrow & H^{n}\left(V_{d_{1}}\right)\\
^{\cong}\downarrow^{\text{·}N} &  & \downarrow^{\text{·}N}\\
H^{n}\left(V_{d_{1}}\cap V_{d_{2}}\right) & \overset{\cong}{\longrightarrow} & H^{n+2}\left(V_{d_{1}}\right).
\end{array}\label{eq:prim}
\end{equation}

The above analysis is identical if the roles of $d_{1}$ and $d_{2}$
are reversed.

\subsection{Accommodating the exceptional locus of the small resolution}

The morphism of exact sequences of mixed Hodge structures
\[
\begin{array}{ccc}
\ldots &  & \ldots\\
\downarrow &  & \downarrow\\
H^{k-2}\left(V_{d_{1}}\cap V_{d_{2}}\right) & \leftrightarrow & H^{k-2}\left(V_{d_{1}}\cap V_{d_{2}}\right)\\
\downarrow &  & \downarrow\\
H^{k}\left(V_{d_{2}}\right) & \rightarrow & H^{k}\left(\tilde{V}_{d_{2}}\right)\\
\downarrow &  & \downarrow\\
H^{k}\left(V_{d_{2}}-\left(V_{d_{1}}\cap V_{d_{2}}\right)\right) & \rightarrow & H^{k}\left(\tilde{V}_{d_{2}}-\left(V_{d_{1}}\cap V_{d_{2}}\right)\right)\\
\downarrow &  & \downarrow\\
H^{k-1}\left(V_{d_{1}}\cap V_{d_{2}}\right) & \leftrightarrow & H^{k-1}\left(V_{d_{1}}\cap V_{d_{2}}\right)\rightarrow\\
\downarrow &  & \downarrow\\
H^{k+1}\left(V_{d_{2}}\right) & \rightarrow & H^{k+1}\left(\tilde{V}_{d_{2}}\right)\\
\downarrow &  & \downarrow\\
\ldots &  & \ldots
\end{array}
\]
together with the isomorphism $H^{k}\left(\tilde{V}_{d_{2}}\right)=H^{k}\left(V_{d_{2}}\right)\oplus H^{k-2}\left(V_{d}\cap V_{d_{1}}\cap V_{d_{2}}\right)\otimes\left(1,1\right)$
yields that for all $k$
\[
H^{k}\left(\tilde{V}_{d_{2}}-\left(V_{d_{1}}\cap V_{d_{2}}\right)\right)=H^{k}\left(V_{d_{2}}-\left(V_{d_{1}}\cap V_{d_{2}}\right)\right)\oplus H^{k-2}\left(V_{d}\cap V_{d_{1}}\cap V_{d_{2}}\right)\otimes\left(1,1\right).
\]

\subsection{Ampleness, the `hard' Lefschetz theorem and the circle bundle $T$}

As we have seen in Subsection \ref{subsec:Mixed-Hodge-structures}
in (\ref{eq:biggest})the right-hand vertical exact sequence for the
circle bundle $T/\left(V_{d_{1}}\cap V_{d_{2}}\right)$ where $N=c_{1}\left(\mathcal{O}_{B}\left(1\right)\right)$
yields
\[
H^{k}\left(T\right)=\left\{ \begin{array}{c}
\mathrm{coker}\left(H^{k-2}\left(V_{d_{1}}\cap V_{d_{2}}\right)\overset{\text{·}d_{1}N}{\longrightarrow}H^{k}\left(V_{d_{1}}\cap V_{d_{2}}\right)\right)+\\
\ker\left(H^{k-1}\left(V_{d_{1}}\cap V_{d_{2}}\right)\overset{\text{·}d_{1}N}{\longrightarrow}H^{k+1}\left(V_{d_{1}}\cap V_{d_{2}}\right)\right).
\end{array}\right.
\]
Since $N$ is ample, it polarizes cohomology so that this formula
reduces to
\begin{equation}
\begin{array}{c}
H^{k}\left(T\right)=H_{prim}^{k}\left(V_{d_{1}}\cap V_{d_{2}}\right)\,\,\,if\,k\leq\dim\left(V_{d_{1}}\cap V_{d_{2}}\right)=n-1\\
H^{k}\left(T\right)=H_{prim}^{k-1}\left(V_{d_{1}}\cap V_{d_{2}}\right)\otimes\left(1,1\right)\,\,\,if\,k\geq\dim\left(V_{d_{1}}\cap V_{d_{2}}\right)+1=n.
\end{array}\label{eq:tee}
\end{equation}
By (\ref{eq:tee}) the long exact cohomology sequence
\[
\begin{array}{c}
\ldots\rightarrow H^{n-1}\left(T\right)\rightarrow H^{n}\left(V_{d}\right)\rightarrow H^{n}\left(V_{d_{1}}-\left(V_{d_{1}}\cap V_{d_{2}}\right)\right)\oplus H^{n}\left(\tilde{V}_{d_{2}}-\left(V_{d_{1}}\cap V_{d_{2}}\right)\right)\\
\rightarrow H^{n}\left(T\right)\rightarrow H^{n+1}\left(V_{d}\right)\rightarrow H^{n+1}\left(V_{d_{1}}-\left(V_{d_{1}}\cap V_{d_{2}}\right)\right)\oplus H^{n+1}\left(\tilde{V}_{d_{2}}-\left(V_{d_{1}}\cap V_{d_{2}}\right)\right)\\
\rightarrow H^{n+1}\left(T\right)\rightarrow H^{n+2}\left(V_{d}\right)\rightarrow\ldots
\end{array}
\]
in (\ref{eq:biggest}) becomes
\begin{equation}
\begin{array}{c}
\ldots\rightarrow H_{prim}^{n-2}\left(V_{d_{1}}\cap V_{d_{2}}\right)\rightarrow H^{n-1}\left(V_{d}\right)\rightarrow H^{n-1}\left(V_{d_{1}}-\left(V_{d_{1}}\cap V_{d_{2}}\right)\right)\oplus H^{n-1}\left(\tilde{V}_{d_{2}}-\left(V_{d_{1}}\cap V_{d_{2}}\right)\right)\\
\rightarrow H_{prim}^{n-1}\left(V_{d_{1}}\cap V_{d_{2}}\right)\rightarrow H^{n}\left(V_{d}\right)\rightarrow H^{n}\left(V_{d_{1}}-\left(V_{d_{1}}\cap V_{d_{2}}\right)\right)\oplus H^{n}\left(\tilde{V}_{d_{2}}-\left(V_{d_{1}}\cap V_{d_{2}}\right)\right)\\
\rightarrow H_{prim}^{n-1}\left(V_{d_{1}}\cap V_{d_{2}}\right)\otimes\left(1,1\right)\rightarrow H^{n+1}\left(V_{d}\right)\rightarrow H^{n+1}\left(V_{d_{1}}-\left(V_{d_{1}}\cap V_{d_{2}}\right)\right)\oplus H^{n+1}\left(\tilde{V}_{d_{2}}-\left(V_{d_{1}}\cap V_{d_{2}}\right)\right)\\
\rightarrow H_{prim}^{n}\left(V_{d_{1}}\cap V_{d_{2}}\right)\otimes\left(1,1\right)\rightarrow H^{n+2}\left(V_{d}\right)\rightarrow\ldots.
\end{array}\label{eq:long}
\end{equation}

\subsection{The monodromy operator\label{subsec:The-monodromy-operator}}

We have the exact sequence of mixed Hodge structures
\begin{equation}
\begin{array}{c}
0\\
\downarrow\\
\mathrm{coker}\left(H^{k}\left(V_{d_{1}}\right)\oplus H^{k}\left(\tilde{V}_{d_{2}}\right)\rightarrow H^{k}\left(V_{d_{1}}\cap V_{d_{2}}\right)\right)\\
\downarrow\\
H^{k+1}\left(V_{d_{1}}\cup\tilde{V}_{d_{2}}\right)\\
\downarrow\\
\mathrm{\ker}\left(H^{k+1}\left(V_{d_{1}}\right)\oplus H^{k+1}\left(\tilde{V}_{d_{2}}\right)\rightarrow H^{k+1}\left(V_{d_{1}}\cap V_{d_{2}}\right)\right)\\
\downarrow\\
0.
\end{array}\label{mixed}
\end{equation}

Consider the commutative diagram 
\[
\begin{array}{ccc}
H^{k}\left(V_{d_{1}}\cap V_{d_{2}}\right) & \rightarrow & H^{k+1}\left(V_{d_{1}}\cup\tilde{V}_{d_{2}}\right)\\
\downarrow^{\rho^{\ast}} &  & \downarrow^{\rho^{\ast}}\\
H^{k}\left(T\right) & \rightarrow & H^{k+1}\left(V_{d}\right)
\end{array}
\]
for $k\geq n$ from (\ref{eq:biggest}). In that case by (\ref{eq:weights2})
the weight of $H^{k}\left(T\right)$ is $k+1$ and so the weight of
its image in $H^{k}\left(V_{d}\right)$ is also $n+1$ implying that
all of the cohomology of $H^{k}\left(V_{d}\right)$ has weight greater
than or equal to $n+1$. Since
\[
H^{k+1}\left(V_{d_{1}}-\left(V_{d_{1}}\cap V_{d_{2}}\right)\right)\oplus H^{k+1}\left(\tilde{V}_{d_{1}}-\left(V_{d_{1}}\cap V_{d_{2}}\right)\right)=H^{k-1}\left(V_{d}\cap V_{d_{1}}\cap V_{d_{2}}\right)\otimes\left(1,1\right)
\]
is of pure weight $k+1$ as is 
\[
H^{k}\left(T\right)=H_{prim}^{k-1}\left(V_{d_{1}}\cap V_{d_{2}}\right)\otimes\left(1,1\right),
\]
$H^{k+1}\left(V_{d}\right)$ must have pure weight $k+1$ for all
$k>n$. However by \cite{Clemens-0}, the unipotent monodromy operator
\begin{equation}
M:H^{k}\left(V_{d}\right)\rightarrow H^{k}\left(V_{d}\right)\label{eq:mon}
\end{equation}
is non-trivial if and only if $\frac{W_{k+1}H^{k}\left(V_{d}\right)}{W_{k}\left(H^{k}\left(V_{d}\right)\right)}$
is non-zero in which case
\begin{equation}
\log M:\frac{W_{k+1}H^{k}\left(V_{d}\right)}{W_{k}\left(H^{k}\left(V_{d}\right)\right)}\longrightarrow\frac{W_{k-1}\left(H^{k}\left(V_{d}\right)\right)}{W_{k-2}\left(H^{k}\left(V_{d}\right)\right)}\label{eq:log-1}
\end{equation}
is an isomorphism. Therefore (\ref{eq:mon}) is the identity map for
$k>n$. Since $M$ preserves the intersection pairing (\ref{eq:mon})
is therefore the identity for all $k\neq n$. Since the mapping $\rho^{\ast}$
in (\ref{eq:biggest}) is the mapping $\mu$, horizontal exactness
at $\mu$ in diagram (3.6) of \cite{Clemens-0} then implies that
the map
\begin{equation}
\rho^{\ast}:H^{k+1}\left(\left(V_{d_{1}}\cup V_{d_{2}}\right)\right)\rightarrow H^{k+1}\left(V_{d}\right)\label{eq:surj}
\end{equation}
is surjective for all $k\neq n$. 

\subsection{Breaking down the governing diagram}

We next state and prove a series of six Lemmas that will establish
the 'building blocks' for the recursive algorithm that is the main
goal of this paper.

From (\ref{eq:tee}), the Lefschetz theorems and the Strictness Lemma
for mixed Hodge structures we conclude the following.
\begin{lem}
\label{lem:i)-In-}In (\ref{eq:biggest}) the mapping
\[
H_{prim}^{k}\left(V_{d_{1}}\cap V_{d_{2}}\right)\overset{\rho^{\ast}}{\longrightarrow}H^{k}\left(T\right)
\]
is an isomorphism if $k\leq n-1$, and the mapping
\[
H^{k}\left(T\right)\overset{\mathrm{Residue}}{\longrightarrow}H_{prim}^{k-1}\left(V_{d_{1}}\cap V_{d_{2}}\right)\otimes\left(1,1\right)
\]
is an isomorphism if $k\geq n$. 
\end{lem}

\begin{proof}
As mentioned above, since $N$ is ample, (\ref{eq:tee}) follows from
the Hard Lefschetz Theorem that implies that the exterior product
with the $\left(1,1\right)$-class $N$ is injective up to the middle
dimension $n-1$ and surjective after that.
\end{proof}
\begin{lem}
\label{lem:ii)-The-mapping} The mapping cone cohomology $\mathbb{H}^{\text{·}+1}\left(\rho^{\ast}\right)$
is isomorphic to $H^{\text{·}-1}\left(V_{d_{1}}\cap V_{d_{2}}\right)\otimes\left(1,1\right)$
in such a way that the horizontal sequence 
\[
\begin{array}{c}
\ldots\rightarrow\mathbb{H}^{\text{·}+1}\left(\rho^{\ast}\right)\rightarrow\\
H^{\text{·}-1}\left(V_{d_{1}}\cap V_{d_{2}}\right)\otimes\left(1,1\right)\oplus H^{\text{·}-1}\left(V_{d_{1}}\cap V_{d_{2}}\right)\otimes\left(1,1\right)\\
\rightarrow H^{\cdot-1}\left(V_{d_{1}}\cap V_{d_{2}}\right)\otimes\left(1,1\right)\rightarrow\ldots
\end{array}
\]
in (\ref{eq:biggest}) in which it sits is the tautological one, that
is, the first map above is the diagonal map and the second is the
subtraction map.
\end{lem}

\begin{proof}
The assertion follows immediately from the fact that $V_{1}\cup\tilde{V}_{2}$
has the same homotopy type as the pasting of $V_{d}$ to the complex
disk bundle over $V_{d_{1}}\cap V_{d_{2}}$ with boundary $T$ to
$T\subseteq V_{d}$ by identifying corresponding points on $T$.
\end{proof}
\begin{lem}
\label{lem:iii)-For-,and}For $k>n$,
\[
\begin{array}{c}
H^{k}\left(V_{d_{1}}-\left(V_{d_{1}}\cap V_{d_{2}}\right)\right)=0\\
H^{k}\left(\tilde{V}_{d_{2}}-\left(V_{d_{1}}\cap V_{d_{2}}\right)\right)\cong H^{k-2}\left(V_{d}\cap V_{d_{1}}\cap V_{d_{2}}\right)\otimes\left(1,1\right)
\end{array}
\]
and the residue mapping
\[
H^{k}\left(\tilde{V}_{d_{2}}-\left(V_{d_{1}}\cap V_{d_{2}}\right)\right)\rightarrow H^{k-1}\left(V_{d_{1}}\cap V_{d_{2}}\right)
\]
is zero. Therefore for $k\geq n+1$ the sequence in (\ref{eq:biggest})
\[
0\rightarrow H^{k}\left(T\right)\rightarrow H^{k+1}\left(V_{d}\right)\rightarrow H^{k-1}\left(V_{d}\cap V_{d_{1}}\cap V_{d_{2}}\right)\otimes\left(1,1\right)\rightarrow0
\]
is exact. Therefore
\[
H^{k+1}\left(V_{d}\right)\cong H^{k+1}\left(V_{1}\right).
\]
\end{lem}

\begin{proof}
Since $\left(V_{d_{1}}\cap V_{d_{2}}\right)$ is ample as a divisor
in $V_{j}$
\[
H^{k}\left(V_{d_{1}}-\left(V_{d_{1}}\cap V_{d_{2}}\right)\right)=0=H^{k}\left(V_{d_{2}}-\left(V_{d_{1}}\cap V_{d_{2}}\right)\right)
\]
for $k>n$. Since the center of the blow-up lies on $\left(V_{d_{1}}\cap V_{d_{2}}\right)$,
$\tilde{V}_{d_{2}}-\left(V_{d_{1}}\cap V_{d_{2}}\right)$ has the
same homotopy type as $V_{d_{2}}-\left(V_{d_{1}}\cap V_{d_{2}}\right)$.
\[
H^{k}\left(\tilde{V}_{d_{2}}-\left(V_{d_{1}}\cap V_{d_{2}}\right)\right)\cong H^{k-2}\left(V_{d}\cap V_{d_{1}}\cap V_{d_{2}}\right)\otimes\left(1,1\right).
\]
For the last assertion, by the Lefschetz Hyperplane Theorem and duality
\[
\begin{array}{c}
H^{k-1}\left(V_{d}\cap V_{d_{1}}\cap V_{d_{2}}\right)\cong H^{2\left(n-2\right)-\left(k-1\right)}\left(V_{d_{1}}\cap V_{d_{2}}\right)\\
\cong H^{2\left(n-2\right)-\left(k-1\right)}\left(V_{d_{2}}\right)\cong H^{k+3}\left(V_{d_{2}}\right)\cong H^{k-3}\left(V_{1}\right)
\end{array}
\]
and
\[
\begin{array}{c}
H^{k}\left(T\right)\cong H_{prim}^{k-1}\left(V_{d_{1}}\cap V_{d_{2}}\right)\cong H_{prim}^{2\left(n-1\right)-\left(k-1\right)}\left(V_{d_{1}}\cap V_{d_{2}}\right)\\
\cong H_{prim}^{2\left(n-1\right)-\left(k-1\right)}\left(V_{d_{2}}\right)\cong H_{prim}^{2\left(n-1\right)-\left(k-1\right)}\left(V_{1}\right)\cong H_{prim}^{k+1}\left(V_{1}\right).
\end{array}
\]
\end{proof}
\begin{lem}
\label{lem:iv)-The-isomorphism}The isomorphism $\frac{H^{n}\left(T\right)}{W_{n}\left(H^{n}\left(T\right)\right)}\rightarrow H_{prim}^{n-1}\left(V_{d_{1}}\cap V_{d_{2}}\right)\otimes\left(1,1\right)$
of Hodge structures of weight $n+1$ induces an additional isomorphism
\[
H_{prim}^{n-1}\left(V_{d_{1}}\cap V_{d_{2}}\right)\otimes\left(1,1\right)\cong\frac{H^{n}\left(T\right)}{W_{n}\left(H^{n}\left(T\right)\right)}\cong\frac{H^{n}\left(V_{d}\right)}{W_{n}\left(H^{n}\left(V_{d}\right)\right)}.
\]
\end{lem}

\begin{proof}
By Subsection \ref{subsec:The-monodromy-operator} $\frac{H^{k}\left(T\right)}{W_{k}\left(H^{k}\left(T\right)\right)}=0$
unless $k=n$. By Lemma \ref{lem:ii)-The-mapping}, the subtraction
map 
\[
H^{\text{·}-1}\left(V_{d_{1}}\cap V_{d_{2}}\right)\otimes\left(1,1\right)\oplus H^{\text{·}-1}\left(V_{d_{1}}\cap V_{d_{2}}\right)\otimes\left(1,1\right)\rightarrow H^{\cdot-1}\left(V_{d_{1}}\cap\tilde{V}_{d_{2}}\right)\otimes\left(1,1\right)
\]
in (\ref{eq:biggest}) is surjective. We first show that $\gamma\in H^{n}\left(T\right)=H_{prim}^{n-1}\left(V_{d_{1}}\cap V_{d_{2}}\right)$
is the difference of residues of some
\begin{equation}
\left(\Gamma_{1},\Gamma_{2}\right)\in H^{n}\left(A_{V_{d_{1}}}^{\cdot}\left(\log\left(V_{d_{1}}\cap V_{d_{2}}\right)\right)\right)\oplus H^{n}\left(A_{\tilde{V}{}_{d_{2}}}^{\text{·}}\left(\log\left(V_{d_{1}}\cap V_{d_{2}}\right)\right)\right).\label{eq:chain}
\end{equation}
By Lemma \ref{lem:iii)-For-,and} we have exact
\[
H^{n}\left(A_{V_{d_{1}}}^{\cdot}\left(\log\left(V_{d_{1}}\cap V_{d_{2}}\right)\right)\right)\rightarrow H^{n-1}\left(V_{d_{1}}\cap V_{d_{2}}\right)\rightarrow H^{n+1}\left(V_{d_{1}}\right)\rightarrow0
\]
and the mapping
\[
H^{n-1}\left(V_{d_{1}}\right)\overset{\text{·}N}{\longrightarrow}H^{n+1}\left(V_{d_{1}}\right)
\]
is a bijection. Also by (\ref{eq:prim})
\[
\frac{H^{n}\left(V_{d_{1}}\right)}{H_{prim}^{n}\left(V_{d_{1}}\right)}
\]
is the kernel of 
\[
H^{n}\left(A_{V_{d_{1}}}^{\cdot}\left(\log\left(V_{d_{1}}\cap V_{d_{2}}\right)\right)\right)\rightarrow H^{n-1}\left(V_{d_{1}}\cap V_{d_{2}}\right)
\]
and by (\ref{eq:push})
\[
\frac{H^{n+1}\left(V_{d_{1}}\right)}{H_{prim}^{n-1}\left(V_{d_{1}}\cap V_{d_{2}}\right)}
\]
is the cokernel of 
\[
H^{n-1}\left(V_{d_{1}}\cap V_{d_{2}}\right)\rightarrow H^{n+1}\left(V_{d_{1}}\right)
\]
Therefore $H_{prim}^{n}\left(V_{d_{1}}\cap V_{d_{2}}\right)$ is the
image of $H^{n}\left(A_{V_{d_{1}}}^{\cdot}\left(\log\left(V_{d_{1}}\cap V_{d_{2}}\right)\right)\right)$
under the residue map. 

The argument for $H^{n}\left(A_{\tilde{V}{}_{d_{2}}}^{\text{·}}\left(\log\left(V_{d_{1}}\cap V_{d_{2}}\right)\right)\right)$
is completely analogous except that $-N$ replaces $N$ reflecting
the reversal of the first Chern class of the normal bundle to $\left(V_{d_{1}}\cap V_{d_{2}}\right)$.

Given non-zero $\gamma\in H^{n}\left(T\right)=H_{prim}^{n-1}\left(V_{d_{1}}\cap V_{d_{2}}\right)$
one can choose $\left(\Gamma_{1},\Gamma_{2}\right)\in H^{k-1}\left(V_{d_{1}}\cap V_{d_{2}}\right)\otimes\left(1,1\right)\oplus H^{\text{k}-1}\left(V_{d_{1}}\cap V_{d_{2}}\right)\otimes\left(1,1\right)$
such that $\mathrm{Residue}\left(\Gamma_{1},\Gamma_{2}\right)$ maps
to $\left(\gamma,\gamma\right)\in H_{prim}^{n-1}\left(V_{d_{1}}\cap V_{d_{2}}\right)$,
it is the image of some non-zero
\begin{equation}
\vartheta\in H^{k}\left(V_{d}\right).\label{eq:cycle}
\end{equation}
The difference of two choices of $\vartheta$ must go to zero in $\mathbb{H}^{k+1}\left(\rho^{\ast}\right)$
since the mapping from $\mathbb{H}^{k+1}\left(\rho^{\ast}\right)$
to $H^{k-1}\left(V_{d_{1}}\cap V_{d_{2}}\right)\otimes\left(1,1\right)\oplus H^{\text{k}-1}\left(V_{d_{1}}\cap V_{d_{2}}\right)\otimes\left(1,1\right)$
is the diagonal mapping. So $\gamma$ determines a well-defined non-zero
element of $\frac{H^{k}\left(V_{d}\right)}{W_{k}\left(H^{k}\left(V_{d}\right)\right)}$.
Surjectivity derives from \cite{Clemens-0} where it is shown that
any non-zero element of $\frac{H^{k}\left(V_{d}\right)}{W_{k}\left(H^{k}\left(V_{d}\right)\right)}$
must have non-trivial local monodromy and therefore give an element
(\ref{eq:chain}) with non-zero but cancelling residues.
\end{proof}
\begin{lem}
\label{lem:v)-The-image}The image of the mapping 
\[
H^{k-1}\left(T\right)\longrightarrow W_{k-1}\left(H^{k}\left(V_{d}\right)\right)
\]
in (\ref{eq:biggest}) lies in $W_{k-1}\left(H^{k}\left(V_{d}\right)\right)$
for all $k\leq n$. This map is an isomorphism 
\[
H_{prim}^{k-1}\left(V_{d_{1}}\cap V_{d_{2}}\right)\longrightarrow W_{k-1}\left(H^{k}\left(V_{d}\right)\right)
\]
for $k=n$. Furthermore
\[
\frac{W_{n}\left(H^{n}\left(V_{d}\right)\right)}{W_{n-1}\left(H^{n}\left(V_{d}\right)\right)}=H_{prim}^{n}\left(V_{d_{1}}\right)\oplus H_{prim}^{n}\left(V_{d_{2}}\right)\oplus H^{n-2}\left(V_{d}\cap V_{d_{1}}\cap V_{d_{2}}\right)\otimes\left(1,1\right)
\]
where
\[
H_{prim}^{n}\left(V_{j}\right)=\ker\left(H^{n}\left(V_{j}\right)\overset{\text{·}N}{\longrightarrow}H^{n+2}\left(V_{j}\right)\right)=\mathrm{coker}\left(H^{n-2}\left(V_{j}\right)\overset{\text{·}N}{\longrightarrow}H^{n}\left(V_{j}\right)\right).
\]
\end{lem}

\begin{proof}
The first assertion follows from i) by weights. By Lemma \ref{lem:iv)-The-isomorphism}
we have that
\[
H_{prim}^{n-1}\left(V_{d_{1}}\cap V_{d_{2}}\right)\otimes\left(1,1\right)\cong\frac{H^{n}\left(V_{d}\right)}{W_{n}\left(H^{n}\left(V_{d}\right)\right)}.
\]
In \cite{Clemens-0} it is shown that weights of classes supported
small close to $V_{d_{1}}\cap V_{d_{2}}$ have weight less than their
degree. Since
\[
H_{prim}^{n-1}\left(V_{d_{1}}\cap V_{d_{2}}\right)\otimes\left(1,1\right)\rightarrow H^{n}\left(V_{d}\right)\rightarrow H^{n}\left(V_{d_{1}}-\left(V_{d_{1}}\cap V_{d_{2}}\right)\right)\oplus H^{n}\left(\tilde{V}_{d_{2}}-\left(V_{d_{1}}\cap V_{d_{2}}\right)\right)
\]
is an exact sequence of mixed Hodge structures and 
\[
W_{k-1}\left(H^{k}\left(V_{d_{1}}-\left(V_{d_{1}}\cap V_{d_{2}}\right)\right)\oplus H^{k}\left(\tilde{V}_{d_{2}}-\left(V_{d_{1}}\cap V_{d_{2}}\right)\right)\right)=0
\]
for all $k,$ $W_{n-1}\left(H^{n}\left(V_{d}\right)\right)=\mathrm{image}\left(H_{prim}^{n-1}\left(V_{d_{1}}\cap V_{d_{2}}\right)\cong H^{n-1}\left(T\right)\rightarrow H^{n}\left(V_{d}\right)\right)$.
Now use the result from \cite{Clemens-0} that (\ref{eq:log-1}) is
an isomorphism. Therefore
\[
\begin{array}{c}
\frac{W_{n}\left(H^{n}\left(V_{d}\right)\right)}{W_{n-1}\left(H^{n}\left(V_{d}\right)\right)}=W_{n}\left(H^{n}\left(V_{d_{1}}-\left(V_{d_{1}}\cap V_{d_{2}}\right)\right)\oplus H^{n}\left(\tilde{V}_{d_{2}}-\left(V_{d_{1}}\cap V_{d_{2}}\right)\right)\right)=\\
W_{n}\left(H^{n}\left(V_{d_{1}}-\left(V_{d_{1}}\cap V_{d_{2}}\right)\right)\oplus H^{n}\left(V_{d_{2}}-\left(V_{d_{1}}\cap V_{d_{2}}\right)\right)\right)\oplus H^{n-2}\left(V_{d}\cap V_{d_{1}}\cap V_{d_{2}}\right)\otimes\left(1,1\right)
\end{array}
\]
and
\[
\begin{array}{c}
W_{n}\left(H^{n}\left(V_{d_{j}}-\left(V_{d_{1}}\cap V_{d_{2}}\right)\right)\right)=\mathrm{image}\left(H^{n}\left(V_{d_{j}}\right)\rightarrow H^{n}\left(V_{d_{j}}-\left(V_{d_{1}}\cap V_{d_{2}}\right)\right)\right)\\
=\frac{H^{n}\left(V_{d_{j}}\right)}{\mathrm{image}\left(H^{n-2}\left(V_{d_{1}}\cap V_{d_{2}}\right)\overset{\text{·}N}{\longrightarrow}H^{n}\left(V_{d_{j}}\right)\right)}=H_{prim}^{n}\left(V_{d_{j}}\right)
\end{array}
\]
by (\ref{eq:prim}).
\end{proof}
\begin{lem}
\label{lem:vi)-If-and}If $p+q=n+1$
\[
H^{n+1}\left(V_{d}\right)\cong\frac{\mathrm{\ker}\left(H^{n+1}\left(V_{d_{1}}\right)\oplus H^{n+1}\left(\tilde{V}_{d_{2}}\right)\rightarrow H^{n+1}\left(V_{d_{1}}\cap V_{d_{2}}\right)\right)}{H^{n+1}\left(V_{d_{1}}\cap V_{d_{2}}\right)}
\]
and is of pure weight $n+1$.
\end{lem}

\begin{proof}
From (\ref{eq:biggest}), Subsection (\ref{subsec:The-monodromy-operator}),
and Lemmas \ref{lem:i)-In-}-\ref{lem:v)-The-image} we have the follow
exact diagram

\[
\begin{array}{ccccc}
\begin{array}{c}
H^{n-1}\left(T\right)\cong\\
H_{prim}^{n-1}\left(V_{d_{1}}\cap V_{d_{2}}\right)
\end{array} & \rightarrow & H^{n}\left(V_{d}\right) & \rightarrow & \begin{array}{c}
H^{n}\left(A_{V_{d_{1}}}^{\text{·}}\left(\log\left(V_{d_{1}}\cap V_{d_{2}}\right)\right)\right)\oplus\\
H^{n}\left(A_{\tilde{V}{}_{d_{2}}}^{\text{·}}\left(\log\left(V_{d_{1}}\cap V_{d_{2}}\right)\right)\right)
\end{array}\\
\downarrow^{\mathrm{Residue}} &  & \downarrow^{\alpha} &  & \downarrow^{\mathrm{Residue}}\\
H^{n-1}\left(V_{d_{1}}\cap V_{d_{2}}\right)\otimes\left(1,1\right) & \overset{0}{\longrightarrow} & H^{n-1}\left(V_{d_{1}}\cap V_{d_{2}}\right)\otimes\left(1,1\right) &  & \begin{array}{c}
H^{n-1}\left(V_{d_{1}}\cap V_{d_{2}}\right)\otimes\left(1,1\right)\oplus\\
H^{n-1}\left(V_{d_{1}}\cap V_{d_{2}}\right)\otimes\left(1,1\right)
\end{array}\\
\downarrow &  & \downarrow\beta &  & \downarrow\\
H^{n}\left(V_{d_{1}}\cap V_{d_{2}}\right) & \overset{\gamma}{\longrightarrow} & H^{n+1}\left(V_{d_{1}}\cup\tilde{V}_{d_{2}}\right) & \rightarrow & \begin{array}{c}
H^{n+1}\left(V_{d_{1}}\right)\oplus\\
H^{n+1}\left(\tilde{V}_{d_{2}}\right)
\end{array}\\
\downarrow^{0} &  & \downarrow^{\rho^{\ast}} &  & \downarrow\\
\begin{array}{c}
H^{n}\left(T\right)\cong\\
H_{prim}^{n-1}\left(V_{d_{1}}\cap V_{d_{2}}\right)\otimes\left(1,1\right)
\end{array} & \rightarrow & H^{n+1}\left(V_{d}\right) & \rightarrow & H^{n-1}\left(V_{d}\cap V_{d_{1}}\cap V_{d_{2}}\right)\otimes\left(1,1\right)\\
 &  & \downarrow &  & 0
\end{array}
\]
where the top right entry maps to the lower left entry to continue
horizontal exactness and 
\begin{equation}
\mathrm{coker\,\alpha}=\mathrm{image\,\beta}=\ker\rho^{\ast}=\frac{H^{n-1}\left(V_{d_{1}}\cap V_{d_{2}}\right)\otimes\left(1,1\right)}{H_{prim}^{n-1}\left(V_{d_{1}}\cap V_{d_{2}}\right)\otimes\left(1,1\right)}\cong H^{n+1}\left(V_{d_{1}}\cap V_{d_{2}}\right)\label{eq:extra}
\end{equation}
of pure weight $n+1$ is the kernel of $\rho^{\ast}$. 
\end{proof}
The sequence of six Lemmas proved just above allow us to prove by
weights and the Strictness Lemma that
\[
\mathrm{coker}\left(H^{n}\left(V_{d_{1}}\right)\oplus H^{n}\left(\tilde{V}_{d_{2}}\right)\rightarrow H^{n}\left(V_{d_{1}}\cap V_{d_{2}}\right)\right)\subseteq\ker\rho^{\ast}.
\]
Therefore the composition
\[
\begin{array}{c}
\mathrm{\ker}\left(H^{n+1}\left(V_{d_{1}}\right)\oplus H^{n+1}\left(\tilde{V}_{d_{2}}\right)\rightarrow H^{n+1}\left(V_{d_{1}}\cap V_{d_{2}}\right)\right)\\
\downarrow\beta\\
H^{n+1}\left(V_{d_{1}}\cup\tilde{V}_{d_{2}}\right)\\
\downarrow^{\rho^{\ast}}\\
H^{n+1}\left(V_{d}\right)
\end{array}
\]
is surjective. So by (\ref{eq:extra})
\[
H^{n+1}\left(V_{d}\right)\cong\frac{\mathrm{\ker}\left(H^{n+1}\left(V_{d_{1}}\right)\oplus H^{n+1}\left(\tilde{V}_{d_{2}}\right)\rightarrow H^{n+1}\left(V_{d_{1}}\cap V_{d_{2}}\right)\right)}{H^{n+1}\left(V_{d_{1}}\cap V_{d_{2}}\right)}
\]
of pure weight $n+1$. 

\subsection{The theorem}

Finally we arrive at the Theorem that is the engine of the recursive
algorithm to computing Hodge numbers of proper complete intersections
as described in Cas One and Case Two at the outset. In summary (\ref{eq:long})
becomes
\begin{equation}
\begin{array}{c}
\ldots\rightarrow H_{prim}^{n-2}\left(V_{d_{1}}\cap V_{d_{2}}\right)\rightarrow H^{n-1}\left(V_{d}\right)\rightarrow H^{n-1}\left(V_{d_{1}}-\left(V_{d_{1}}\cap V_{d_{2}}\right)\right)\oplus H^{n-1}\left(\tilde{V}_{d_{2}}-\left(V_{d_{1}}\cap V_{d_{2}}\right)\right)\\
\rightarrow H_{prim}^{n-1}\left(V_{d_{1}}\cap V_{d_{2}}\right)\rightarrow H^{n}\left(V_{d}\right)\rightarrow H^{n}\left(V_{d_{1}}-\left(V_{d_{1}}\cap V_{d_{2}}\right)\right)\oplus H^{n}\left(\tilde{V}_{d_{2}}-\left(V_{d_{1}}\cap V_{d_{2}}\right)\right)\\
\rightarrow H_{prim}^{n-1}\left(V_{d_{1}}\cap V_{d_{2}}\right)\otimes\left(1,1\right)\rightarrow H^{n+1}\left(V_{d}\right)\rightarrow H^{n-1}\left(V_{d}\cap V_{d_{1}}\cap V_{d_{2}}\right)\otimes\left(1,1\right)\rightarrow0
\end{array}\label{eq:longer}
\end{equation}
together with the short exact sequences
\begin{equation}
0\rightarrow H_{prim}^{k-2}\left(V_{d_{1}}\cap V_{d_{2}}\right)\otimes\left(1,1\right)\rightarrow H^{k}\left(V_{d}\right)\rightarrow H^{k-2}\left(V_{d}\cap V_{d_{1}}\cap V_{d_{2}}\right)\otimes\left(1,1\right)\rightarrow0\label{eq:least}
\end{equation}
for $k>n+1$.

By \ref{lem:i)-In-}vi) if $p+q=n+1$
\[
H^{n+1}\left(V_{d}\right)\cong\frac{\mathrm{\ker}\left(H^{n+1}\left(V_{d_{1}}\right)\oplus H^{n+1}\left(\tilde{V}_{d_{2}}\right)\rightarrow H^{n+1}\left(V_{d_{1}}\cap V_{d_{2}}\right)\right)}{H^{n+1}\left(V_{d_{1}}\cap V_{d_{2}}\right)}.
\]

For $k=n$ we have by \ref{lem:i)-In-}v) that 
\begin{equation}
\begin{array}{c}
W_{n-1}\left(H^{n}\left(V_{d}\right)\right)=H_{prim}^{n-1}\left(V_{d_{1}}\cap V_{d_{2}}\right)\\
\frac{W_{n}\left(H^{n}\left(V_{d}\right)\right)}{W_{n-1}\left(H^{n}\left(V_{d}\right)\right)}=H_{prim}^{n}\left(V_{d_{1}}\right)\oplus H_{prim}^{n}\left(V_{d_{2}}\right)\oplus H^{n-2}\left(V_{d}\cap V_{d_{1}}\cap V_{d_{2}}\right)\otimes\left(1,1\right)\\
\frac{W_{n+1}\left(H^{n}\left(V_{d}\right)\right)}{W_{n}\left(H^{n}\left(V_{d}\right)\right)}=H_{prim}^{n-1}\left(V_{d_{1}}\cap V_{d_{2}}\right)\otimes\left(1,1\right)
\end{array}\label{eq:middle}
\end{equation}
gives precisely the full weight filtration of the mixed Hodge structure
on $H^{n}\left(V_{d}\right)$.

Since the dimensions of the associated graded for the Hodge filtrations
on the asymptotic mixed Hodge structure on $H^{\cdot}\left(V_{d}\right)$
are the same as those for the usual Hodge filtrations on $H^{\cdot}\left(V_{d}\right)$
, we have the following result which is somewhat complicated to state
but, as the example in the Introduction shows, very easy to use.
\begin{thm}
\label{thm:The-Hodge-numbers} Case One: For a very ample line bundle
$L^{-1}$ over a projective manifold $B$ of complex dimension $n$
we consider sections $H^{0}\left(L^{-d}\right)$ as proper subvarieties
$V_{d}$ of $\left|L\right|$ via the $d$-th power mapping 
\[
\left|L\right|\rightarrow\left|L^{d}\right|.
\]
Thus $V_{1}\cong B$.

In Case Two, $Y$ is any smooth complex projective manifold and
\[
V_{1}=Y\cap H\,\,\,\dim Y=n+1
\]
is a smooth hyperplane section. 

In both cases let $F_{d}$ denote the homogeneous form defining $V_{d}$.
Then the cohomology groups and the Hodge numbers of a (generic) complete
intersection
\[
V_{d_{1}}\cap V_{d_{2}}\cap\ldots\cap V_{d_{r}}=:V_{d_{1}}\text{·}V_{d_{2}}\text{·}\ldots\text{·}V_{d_{r}}
\]
with $1\leq r\leq n$ are given inductively in terms of the Hodge
numbers of the sequence $V_{1},V_{1}^{2},\ldots,V_{1}^{n}$ by the
following formulas for $d=d_{1}+d_{2}$ where $\left(1,1\right)$
denotes the standard Hodge-Tate structure of dimension one and weight
two:

\noindent i) For cohomology classes in degree $k>n+1$,
\[
H^{k}\left(V_{d}\right)\cong H^{k}\left(V_{1}\right),
\]
of pure weight $k$, 

\noindent for cohomology classes of degree $n+1$
\[
H^{n+1}\left(V_{d}\right)\cong\frac{\ker\left(\left(H^{n+1}\left(V_{d_{1}}\right)\oplus H^{n+1}\left(\tilde{V}_{d_{2}}\right)\right)\rightarrow H^{n+1}\left(V_{d_{1}}\cap V_{d_{2}}\right)\right)}{H^{n+1}\left(V_{d_{1}}\cap V_{d_{2}}\right)},
\]

\noindent and for cohomology classes of degree $n$
\[
H^{n}\left(V_{d}\right)=\left\{ \begin{array}{c}
\frac{H^{n}\left(V_{d}\right)}{W_{n}\left(H^{n}\left(V_{d}\right)\right)}=H_{prim}^{n-1}\left(V_{d_{1}}\cap V_{d_{2}}\right)\otimes\left(1,1\right)\\
\frac{W_{n}\left(H^{n}\left(V_{d}\right)\right)}{W_{n-1}\left(H^{n}\left(V_{d}\right)\right)}=H_{prim}^{n}\left(V_{1}\right)\oplus H_{prim}^{n}\left(V_{2}\right)\oplus H^{n-2}\left(V_{d}\cap V_{d_{1}}\cap V_{d_{2}}\right)\otimes\left(1,1\right)\\
W_{n-1}\left(H^{n}\left(V_{d}\right)\right)=H_{prim}^{n-1}\left(V_{d_{1}}\cap V_{d_{2}}\right).
\end{array}\right.
\]

\noindent ii) If $p+q=k>n+1$
\[
h^{p,q}\left(V_{d}\right)=h^{p,q}\left(V_{1}\right),
\]
if $k=p+q=n+1$
\[
\begin{array}{c}
h^{p,q}\left(V_{d}\right)=h^{p,q}\left(\ker\left(H^{n+1}\left(V_{d_{1}}\right)\oplus H^{n+1}\left(\tilde{V}_{d_{2}}\right)\rightarrow H^{n+1}\left(V_{d_{1}}\cap V_{d_{2}}\right)\right)\right)-h^{p,q}\left(V_{d_{1}}\cap V_{d_{2}}\right)\\
=\dim\left(\ker\left(H^{p,q}\left(V_{d_{1}}\right)\oplus H^{p,q}\left(V_{d_{2}}\right)\rightarrow H^{p,q}\left(V_{d_{1}}\cap V_{d_{2}}\right)\right)\right)+h^{p-1,q-1}\left(V_{d}\cap V_{d_{1}}\cap V_{d_{2}}\right)-h^{p,q}\left(V_{d_{1}}\cap V_{d_{2}}\right),
\end{array}
\]

\noindent and if $p+q=k=n$,
\[
h^{p,q}\left(V_{d}\right)=\left\{ \begin{array}{c}
h_{prim}^{p-1,q}\left(V_{d_{1}}\cap V_{d_{2}}\right)\\
h_{prim}^{p,q}\left(V_{d_{1}}\right)+h_{prim}^{p,q}\left(V_{d_{2}}\right)+h^{p-1,q-1}\left(V_{d}\cap V_{d_{1}}\cap V_{d_{2}}\right)\\
+h_{prim}^{p,q-1}\left(V_{d_{1}}\cap V_{d_{2}}\right)
\end{array}\right.
\]
where the row of each summand indicates weight in the asymptotic mixed
Hodge structure.
\end{thm}

\appendix

\section{\label{subsec:Computing-Hodge-numbers}Hodge numbers of proper complete
intersections in $Y=\left|L\right|$ for a very ample line bundle
$L^{-1}$ on a complex projective threefold $B$ }

We begin with the linear table notation:

\smallskip{}

\begin{center}
\begin{tabular}{|c|c|c|c|c|}
\hline 
 & $V_{1}^{4}$ & $V_{1}^{3}$ & $V_{1}^{2}$ & $V_{1}=B$\tabularnewline
\hline 
\hline 
$h^{0,0}$ & $l_{0}^{0,0}$ & $l_{1}^{0,0}$ & $l_{2}^{0,0}$ & $l_{3}^{0,0}$\tabularnewline
\hline 
$h^{1,0}$ &  & $l_{1}^{1,0}$ & $l_{2}^{1,0}$ & $l_{3}^{1,0}$\tabularnewline
\hline 
$h^{2,0}$ &  &  & $l_{2}^{2,0}$ & $l_{3}^{2,0}$\tabularnewline
\hline 
$h^{1,1}$ &  &  & $l_{2}^{1,1}$ & $l_{3}^{1,1}$\tabularnewline
\hline 
$h^{3,0}$ &  &  &  & $l_{3}^{3,0}$\tabularnewline
\hline 
$h^{2,1}$ &  &  &  & $l_{3}^{2,1}$\tabularnewline
\hline 
\end{tabular}
\par\end{center}

\smallskip{}

\subsection{Hodge numbers of curves $\left(V_{d_{0}}\cap V_{d_{1}}\cap V_{d_{2}}\right)\subseteq\left|L\right|\subseteq Q$}

As above let $N=c_{1}\left(L^{-1}\right)$. Now the canonical bundle
of $\left|L\right|$ is 
\[
p^{\ast}\left(K_{B}\otimes L^{-1}\right)
\]
so by adjunction
\[
K_{B}\cong K_{V_{1}}=p^{\ast}\left(K_{B}\otimes L^{-1}\right)\otimes\mathcal{N}_{V_{1}|\left|L\right|}.
\]
Therefore the normal bundle satisfies
\[
\begin{array}{c}
\mathcal{N}_{V_{1}|\left|L\right|}=\left.p^{\ast}L\right|_{V_{1}}=\mathcal{O}_{V_{1}}\otimes p^{\ast}\mathcal{O}_{B}\left(-N\right)\\
\mathcal{N}_{\left(V_{1}\cap V_{1}'\right)|V_{1}}=\mathcal{O}_{V_{1}\cap V_{1}'}\otimes p^{\ast}\mathcal{O}_{B}\left(N\right)=\mathcal{O}_{V_{1}\cap V_{1}'}\left(N\right)\\
\mathcal{N}_{\left(V_{1}\cap V_{1}'\cap V''_{1}\right)|\left(V_{1}\cap V_{1}'\right)}=\mathcal{O}_{V_{1}\cap V_{1}'\cap V''_{1}}\left(2N\right)\\
\left|V_{1}\cap V_{1}'\cap V''_{1}\cap V'''_{1}\right|=N^{3}=:l_{0}^{0,0}
\end{array}
\]
and 
\[
\begin{array}{c}
K_{V_{1}}=p^{\ast}K_{B}\otimes\mathcal{O}_{V_{1}}\\
K_{V_{1}\cap V_{1}'}=p^{\ast}K_{B}\otimes\mathcal{O}_{V_{1}\cap V_{1}'}\left(N\right)\\
K_{V_{1}\cap V_{1}'\cap V''_{1}}=p^{\ast}K_{B}\otimes\mathcal{O}_{V_{1}\cap V_{1}'\cap V''_{1}}\left(2N\right).
\end{array}
\]

Then, abusing notation slightly,
\[
\begin{array}{c}
h^{0}\left(K_{V_{1}\cap V_{1}\cap V_{d}}\right)=h^{0}\left(K_{V_{1}\cap V_{1}\cap V_{d-1}}\right)+h^{0}\left(K_{V_{1}\cap V_{1}\cap V_{1}}\right)+\left(d-1\right)l_{0}^{0,0}-1\\
=d\text{·}h^{0}\left(K_{V_{1}\cap V_{1}\cap V_{1}}\right)+\left(\sum_{k=1}^{d-1}k\right)l_{0}^{0,0}-\left(d-1\right)\\
=d\text{·}l_{1}^{1,0}+\left(\frac{d^{2}-d}{2}\right)\text{·}l_{0}^{0,0}-\left(d-1\right).
\end{array}
\]
Similarly
\[
\begin{array}{c}
h^{0}\left(K_{V_{1}\cap V_{d_{1}}\cap V_{d_{2}}}\right)=h^{0}\left(K_{V_{1}\cap V_{1}\cap V_{d_{2}}}\right)+h^{0}\left(K_{V_{1}\cap V_{d_{1}-1}\cap V_{d_{2}}}\right)+\left(d_{1}-1\right)d_{2}\text{·}l_{0}^{0,0}-1\\
=d_{1}h^{0}\left(K_{V_{1}\cap V_{1}\cap V_{d_{2}}}\right)+\left(\sum_{k=1}^{d_{1}-1}k\right)d_{2}\text{·}l_{0}^{0,0}-\left(d_{1}-1\right)\\
=d_{1}\left(d_{2}\text{·}l_{1}^{1,0}+\left(\frac{d_{2}^{2}-d_{2}}{2}\right)\text{·}l_{0}^{0,0}-\left(d_{2}-1\right)\right)+\left(\frac{d_{1}^{2}-d_{1}}{2}\right)d_{2}\text{·}l_{0}^{0,0}-\left(d_{1}-1\right)
\end{array}
\]
and finally
\[
\begin{array}{c}
h^{0}\left(K_{V_{d_{0}}\cap V_{d_{1}}\cap V_{d_{2}}}\right)=h^{0}\left(K_{V_{1}\cap V_{d_{1}}\cap V_{d_{2}}}\right)+h^{0}\left(K_{V_{d_{0}-1}\cap V_{d_{1}}\cap V_{d_{2}}}\right)+\left(d_{0}-1\right)d_{1}d_{2}l_{0}^{0,0}-1\\
=d_{0}h^{0}\left(K_{V_{1}\cap V_{d_{1}}\cap V_{d_{2}}}\right)+\left(\sum_{k=1}^{d_{0}-1}k\right)d_{1}d_{2}\text{·}l_{0}^{0,0}-\left(d_{0}-1\right)\\
=d_{0}\left(d_{1}\left(d_{2}\text{·}l_{1}^{1,0}+\left(\frac{d_{2}^{2}-d_{2}}{2}\right)\text{·}l_{0}^{0,0}-\left(d_{2}-1\right)\right)+\left(\frac{d_{1}^{2}-d_{1}}{2}\right)d_{2}\text{·}l_{0}^{0,0}-\left(d_{1}-1\right)\right)\\
+\left(\frac{d_{0}^{2}-d_{0}}{2}\right)d_{1}d_{2}\text{·}l_{0}^{0,0}-\left(d_{0}-1\right)
\end{array}
\]
and so
\begin{equation}
h^{1,0}\left(V_{d_{0}}\text{·}V_{d_{1}}\text{·}V_{d_{2}}\right)=d_{0}d_{1}d_{2}\text{·}l_{1}^{1,0}+\left(\frac{d_{0}d_{1}d_{2}\left(d_{0}+d_{1}+d_{2}-3\right)}{2}\right)\text{·}l_{0}^{0,0}-d_{0}d_{1}d_{2}+1\label{eq:curvef}
\end{equation}

In summary we have the following table of point set and curve Hodge
numbers:

\smallskip{}

\begin{center}
\begin{tabular}{|c|c|c|}
\hline 
 & $V_{d_{0}}\text{·}V_{d_{1}}\text{·}V_{d_{2}}\text{·}V_{d_{3}}$ & $V_{d_{0}}\cap V_{d_{1}}\cap V_{d_{2}}\subseteq\left|L\right|\subseteq Q$\tabularnewline
\hline 
\hline 
$h^{0,0}$ & $l_{0}^{0,0}d_{0}\text{·}d_{1}\text{·}d_{2}\text{·}d_{3}$ & $1$\tabularnewline
\hline 
$h^{1,0}$ &  & $d_{0}d_{1}d_{2}\text{·}l_{1}^{1,0}+\left(\frac{d_{0}d_{1}d_{2}\left(d_{0}+d_{1}+d_{2}-3\right)}{2}\right)\text{·}l_{0}^{0,0}-d_{0}d_{1}d_{2}+1$\tabularnewline
\hline 
\end{tabular}
\par\end{center}

\smallskip{}

\subsection{Hodge numbers of surfaces $\left(V_{d_{0}}\cap V_{d_{2}}\right)\subseteq\left|L\right|\subseteq Q$}

To compute the Hodge numbers of surfaces $V_{d_{0}}\cap V_{d}$ we
study the linear degeneration to $V_{d_{0}}\cap\left(V_{d_{1}}\cup V_{d_{2}}\right)$
where $d_{1}+d_{2}=d$. This time we apply Case Two of Theorem \ref{thm:In-either-Case}
with $Y=V_{d_{0}}$ . The total space of the family
\[
\left\{ t\text{·}F_{d}-F_{d_{1}}\text{·}F_{d_{2}}=\left|\begin{array}{cc}
t & F_{d_{1}}\\
F_{d_{2}} & F_{d}
\end{array}\right|=0\right\} \subseteq\mathbb{C}\times V_{d_{0}}
\]
has a nodal locus at the $l_{0}^{0,0}\text{·}d_{0}\text{·}d_{1}\text{·}d_{2}\text{·}d$
points
\[
V_{d_{0}}\cap V_{d_{1}}\cap V_{d_{2}}\cap V_{d}
\]
where all four entries in the above matrix are zero. 

Case Two of Theorem \ref{thm:In-either-Case}with $Y=V_{d_{0}}$ and
the fact that all divisors are very ample then give that for $k=3$
and $Y=V_{d_{0}}$ we have the surface formula
\[
\begin{array}{c}
h^{2,1}\left(V_{d_{0}}\cap V_{d}\right)=\\
\dim\ker\left(H^{2,1}\left(V_{d_{0}}\cap V_{d_{1}}\right)\oplus H^{2,1}\left(V_{d_{0}}\cap\tilde{V}_{d_{2}}\right)\rightarrow H^{2,1}\left(V_{d_{0}}\cap V_{d_{1}}\cap V_{d_{2}}\right)\right)-h^{2,1}\left(V_{d_{0}}\cap V_{d_{1}}\cap V_{d_{2}}\right)\\
=\dim\left(\ker\left(H^{2,1}\left(V_{d_{0}}\cap V_{d_{1}}\right)\oplus H^{2,1}\left(V_{d_{0}}\cap V_{d_{2}}\right)\rightarrow H^{2,1}\left(V_{d_{0}}\cap V_{d_{1}}\cap V_{d_{2}}\right)\right)\right)\\
+h^{1,0}\left(V_{d_{0}}\cap V_{d}\cap V_{d_{1}}\cap V_{d_{2}}\right)-h^{2,1}\left(V_{d_{0}}\cap V_{d_{1}}\cap V_{d_{2}}\right).
\end{array}
\]
Therefore by dimension
\[
h^{2,1}\left(V_{d_{0}}\cap V_{d}\right)=h^{2,1}\left(V_{d_{0}}\cap V_{d_{1}}\right)+h^{2,1}\left(V_{d_{0}}\cap V_{d_{2}}\right)
\]
and so by induction for $Y=V_{d_{0}}$ 
\[
h^{2,1}\left(V_{d_{0}}\cap V_{d_{1}}\right)=d_{1}\text{·}h^{2,1}\left(V_{d_{0}}\cap V_{1}\right)=d_{0}\text{·}d_{1}\text{·}l_{2}^{1,0}.
\]
since $h^{2,1}\left(V_{1}\cap V_{1}\right)=l_{2}^{1,0}$.
\begin{flushleft}
When $p+q=2$, we can use Case Two of Theorem \ref{thm:In-either-Case}
$Y=V_{d_{0}}$ since the formulae
\begin{equation}
\begin{array}{c}
h^{2,0}\left(V_{d_{0}}\cap V_{d}\right)=\left\{ \begin{array}{c}
h^{1,0}\left(V_{d_{0}}\cap V_{d_{1}}\cap V_{d_{2}}\right)+\\
h^{2,0}\left(V_{d_{0}}\cap V_{d_{1}}\right)+h^{2,0}\left(V_{d_{0}}\cap V_{d_{2}}\right)
\end{array}\right\} \\
h^{1,1}\left(V_{d_{0}}\cap V_{d}\right)=\left\{ \begin{array}{c}
h^{0,1}\left(V_{d_{0}}\cap V_{d_{1}}\cap V_{d_{2}}\right)+\\
h_{prim}^{1,1}\left(V_{d_{0}}\cap V_{d_{1}}\right)+h_{prim}^{1,1}\left(V_{d_{0}}\cap V_{d_{2}}\right)+h^{0,0}\left(V_{d_{0}}\cap V_{d_{1}}\cap V_{d_{2}}\cap V_{d}\right)\\
+h^{1,0}\left(V_{d_{0}}\cap V_{d_{1}}\cap V_{d_{2}}\right)
\end{array}\right\} 
\end{array}\label{eq:family}
\end{equation}
allows for a complete induction starting with $h^{p,q}\left(V_{1}\cap V_{1}\right)=l_{2}^{p,q}$.
In particular
\[
\begin{array}{c}
h^{2,0}\left(V_{1}\cap V_{d}\right)=\left\{ \begin{array}{c}
h^{1,0}\left(V_{1}\cap V_{1}\cap V_{d-1}\right)+\\
h^{2,0}\left(V_{1}\cap V_{1}\right)+h^{2,0}\left(V_{1}\cap V_{d-1}\right)
\end{array}\right\} \\
\\
\end{array}.
\]
Therefore if $d>1$
\[
\begin{array}{c}
h^{2,0}\left(V_{1}\cap V_{d}\right)=\left\{ \begin{array}{c}
\left(d-1\right)\text{·}l_{1}^{1,0}+\left(\frac{\left(d-1\right)^{2}-\left(d-1\right)}{2}\right)\text{·}l_{0}^{0,0}-\left(d-1\right)+1\\
+l_{2}^{2,0}+h^{2,0}\left(V_{1}\cap V_{d-1}\right)
\end{array}\right\} \\
h^{2,0}\left(V_{1}\cap V_{d}\right)=\left\{ \begin{array}{c}
d\text{·}l_{2}^{2,0}+\left(\frac{d^{2}-d}{2}\right)l_{1}^{1,0}\\
+\left(\frac{\frac{2\left(d-1\right)^{3}+3\left(d-1\right)^{2}+\left(d-1\right)}{6}-\frac{d^{2}-d}{2}}{2}\right)\text{·}l_{0}^{0,0}-\left(\frac{d^{2}-d}{2}\right)+\left(d-1\right)
\end{array}\right\} \\
h^{2,0}\left(V_{1}\cap V_{d}\right)=\left\{ \begin{array}{c}
d\text{·}l_{2}^{2,0}+\left(\frac{d^{2}-d}{2}\right)l_{1}^{1,0}\\
+\left(\frac{\frac{2\left(d-1\right)^{3}+3\left(d-1\right)^{2}+\left(d-1\right)}{6}-\frac{d^{2}-d}{2}}{2}\right)\text{·}l_{0}^{0,0}\\
-\frac{d^{2}-3d+2}{2}
\end{array}\right\} \\
h^{2,0}\left(V_{1}\cap V_{d}\right)=\left\{ \begin{array}{c}
d\text{·}l_{2}^{2,0}+\left(\frac{d^{2}-d}{2}\right)l_{1}^{1,0}\\
+\left(\frac{d^{3}-3d^{2}+2d}{6}\right)\text{·}l_{0}^{0,0}\\
-\frac{d^{2}-3d+2}{2}
\end{array}\right\} .
\end{array}
\]
Also
\[
h^{1,1}\left(V_{1}\cap V_{d}\right)=\left\{ \begin{array}{c}
h^{0,1}\left(V_{1}\cap V_{1}\cap V_{d-1}\right)+\\
h_{prim}^{1,1}\left(V_{1}\cap V_{1}\right)+h_{prim}^{1,1}\left(V_{1}\cap V_{d-1}\right)+h^{0,0}\left(V_{1}\cap V_{1}\cap V_{d-1}\cap V_{d}\right)\\
+h^{1,0}\left(V_{1}\cap V_{1}\cap V_{d-1}\right)
\end{array}\right\} 
\]
\[
\begin{array}{c}
h^{1,1}\left(V_{1}\cap V_{d}\right)=\left\{ \begin{array}{c}
2\left(\left(d-1\right)\text{·}l_{1}^{1,0}+\left(\frac{\left(d-1\right)\left(d-2\right)}{2}\right)\text{·}l_{0}^{0,0}-\left(d-1\right)+1\right)\\
\left(l_{2}^{1,1}-1\right)+\left(h_{prim}^{1,1}\left(V_{1}\cap V_{d-1}\right)\right)+d\left(d-1\right)\text{·}l_{0}^{0,0}
\end{array}\right\} \\
h_{prim}^{1,1}\left(V_{1}\cap V_{d}\right)=\left\{ \begin{array}{c}
\left(d-1\right)\text{·}l_{1}^{1,0}+\left(\left(d-1\right)^{2}-\left(d-1\right)\right)\text{·}l_{0}^{0,0}-2\left(d-1\right)+1\\
\left(l_{2}^{1,1}-1\right)+\left(h_{prim}^{1,1}\left(V_{1}\cap V_{d-1}\right)\right)+\left(\left(d-1\right)^{2}+\left(d-1\right)\right)\text{·}l_{0}^{0,0}
\end{array}\right\} \\
h_{prim}^{1,1}\left(V_{1}\cap V_{d}\right)=\left\{ \begin{array}{c}
\left(d-1\right)\text{·}l_{1}^{1,0}+2\left(d-1\right)^{2}\text{·}l_{0}^{0,0}-2\left(d-1\right)+1\\
\left(l_{2}^{1,1}-1\right)+\left(h_{prim}^{1,1}\left(V_{1}\cap V_{d-1}\right)\right)
\end{array}\right\} \\
h_{prim}^{1,1}\left(V_{1}\cap V_{d}\right)=\left\{ \begin{array}{c}
d\text{·}\left(l_{2}^{1,1}-1\right)+2\left(\sum_{k=1}^{d-1}k\right)\text{·}l_{1}^{1,0}\\
+2\left(\sum_{k=1}^{d-1}k^{2}\right)\text{·}l_{0}^{0,0}\\
-2\left(\sum_{k=1}^{d-1}k\right)+\left(d-1\right)
\end{array}\right\} \\
h_{prim}^{1,1}\left(V_{1}\cap V_{d}\right)=\left\{ \begin{array}{c}
d\text{·}\left(l_{2}^{1,1}-1\right)+\left(d^{2}-d\right)\text{·}l_{1}^{1,0}\\
+\left(\frac{2d^{3}-3d^{2}+d}{3}\right)\text{·}l_{0}^{0,0}\\
-\left(d^{2}-d\right)+d-1
\end{array}\right\} \\
h^{1,1}\left(V_{1}\cap V_{d}\right)=\left\{ \begin{array}{c}
d\text{·}\left(l_{2}^{1,1}-1\right)+\left(d^{2}-d\right)l_{1}^{1,0}\\
+\left(\frac{2d^{3}-3d^{2}+d}{3}\right)\text{·}l_{0}^{0,0}\\
-d^{2}+2d
\end{array}\right\} .
\end{array}
\]
\par\end{flushleft}

So finally for all positive $d$
\[
\begin{array}{c}
h^{2,0}\left(V_{1}\cap V_{d}\right)=d\text{·}l_{2}^{2,0}+\left(\frac{d^{2}-d}{2}\right)\text{·}l_{1}^{1,0}+\left(\frac{d^{3}-3d^{2}+2d}{6}\right)\text{·}l_{0}^{0,0}-\frac{d^{2}-3d+2}{2}\\
h^{1,1}\left(V_{1}\cap V_{d}\right)=d\text{·}l_{2}^{1,1}+\left(d^{2}-d\right)l_{1}^{1,0}+\left(\frac{2d^{3}-3d^{2}+d}{3}\right)\text{·}l_{0}^{0,0}-d^{2}+d.
\end{array}
\]

Lastly we reverse the roles of $d_{0}$ and $d$ and write
\[
h^{2,0}\left(V_{d_{0}}\cap V_{d}\right)=\left\{ \begin{array}{c}
h^{1,0}\left(V_{d_{0}}\cap V_{d-1}\cap V_{1}\right)+\\
h^{2,0}\left(V_{d_{0}}\cap V_{d-1}\right)+h^{2,0}\left(V_{d_{0}}\cap V_{1}\right)
\end{array}\right\} 
\]
and so for $d>1$
\[
h^{2,0}\left(V_{d_{0}}\cap V_{d}\right)=\left\{ \begin{array}{c}
d_{0}\left(d-1\right)\text{·}l_{1}^{1,0}+d_{0}\left(\frac{\left(d-1\right)\left(d_{0}+d-3\right)}{2}\right)\text{·}l_{0}^{0,0}-d_{0}\left(d-1\right)+1+\\
h^{2,0}\left(V_{d_{0}}\cap V_{d-1}\right)+\left(d_{0}\text{·}l_{2}^{2,0}+\left(\frac{d_{0}^{2}-d_{0}}{2}\right)\text{·}l_{1}^{1,0}+\left(\frac{d_{0}^{3}-3d_{0}^{2}+2d_{0}}{6}\right)\text{·}l_{0}^{0,0}-\frac{d_{0}^{2}-3d_{0}+2}{2}\right)
\end{array}\right\} 
\]
\[
h^{2,0}\left(V_{d_{0}}\cap V_{d}\right)=\left\{ \begin{array}{c}
h^{2,0}\left(V_{d_{0}}\cap V_{d-1}\right)+\left(d_{0}\text{·}l_{2}^{2,0}+\left(\frac{d_{0}^{2}-d_{0}}{2}\right)\text{·}l_{1}^{1,0}+\left(\frac{d_{0}^{3}-3d_{0}^{2}+2d_{0}}{6}\right)\text{·}l_{0}^{0,0}-\frac{d_{0}^{2}-3d_{0}+2}{2}\right)\\
+d_{0}\left(d-1\right)\text{·}l_{1}^{1,0}+d_{0}\left(\frac{\left(d-1\right)^{2}+\left(d_{0}-2\right)\left(d-1\right)}{2}\right)\text{·}l_{0}^{0,0}\\
-d_{0}\left(d-1\right)+1
\end{array}\right\} 
\]
\[
h^{2,0}\left(V_{d_{0}}\cap V_{d}\right)=\left\{ \begin{array}{c}
d\left(d_{0}\text{·}l_{2}^{2,0}+\left(\frac{d_{0}^{2}-d_{0}}{2}\right)\text{·}l_{1}^{1,0}+\left(\frac{d_{0}^{3}-3d_{0}^{2}+2d_{0}}{6}\right)\text{·}l_{0}^{0,0}-\frac{d_{0}^{2}-3d_{0}+2}{2}\right)\\
+d_{0}\left(\sum_{k=1}^{d-1}k\right)\text{·}l_{1}^{1,0}+d_{0}\left(\frac{\left(\sum_{k=1}^{d-1}k^{2}\right)+\left(d_{0}-2\right)\left(\sum_{k=1}^{d-1}k\right)}{2}\right)\text{·}l_{0}^{0,0}\\
-d_{0}\left(\sum_{k=1}^{d-1}k\right)+\left(d-1\right)
\end{array}\right\} 
\]
\[
h^{2,0}\left(V_{d_{0}}\cap V_{d}\right)=\left\{ \begin{array}{c}
d\left(d_{0}\text{·}l_{2}^{2,0}+\left(\frac{d_{0}^{2}-d_{0}}{2}\right)\text{·}l_{1}^{1,0}+\left(\frac{d_{0}^{3}-3d_{0}^{2}+2d_{0}}{6}\right)\text{·}l_{0}^{0,0}-\frac{d_{0}^{2}-3d_{0}+2}{2}\right)\\
+d_{0}\left(\frac{d^{2}-d}{2}\right)\text{·}l_{1}^{1,0}+d_{0}\left(\frac{\left(\frac{2d^{3}-3d^{2}+d}{6}\right)+\left(d_{0}-2\right)\left(\frac{d^{2}-d}{2}\right)}{2}\right)\text{·}l_{0}^{0,0}\\
-d_{0}\left(\frac{d^{2}-d}{2}\right)+\left(d-1\right)
\end{array}\right\} 
\]
\[
h^{2,0}\left(V_{d_{0}}\cap V_{d}\right)=\left\{ \begin{array}{c}
dd_{0}\text{·}l_{2}^{2,0}+\left(d\left(\frac{d_{0}^{2}-d_{0}}{2}\right)+d_{0}\left(\frac{d^{2}-d}{2}\right)\right)\text{·}l_{1}^{1,0}\\
+\left(\frac{2dd_{0}^{3}-6dd_{0}^{2}+4d_{0}d}{12}+\frac{2d_{0}d^{3}-3d_{0}d^{2}+d_{0}d}{12}+\left(\frac{3d_{0}^{2}d^{2}-3d_{0}^{2}d}{12}\right)-\left(\frac{6d_{0}d^{2}-6d_{0}d}{12}\right)\right)\text{·}l_{0}^{0,0}\\
-\left(\frac{d_{0}^{2}d-3d_{0}d+2d}{2}+\frac{d_{0}d^{2}-d_{0}d}{2}\right)+\frac{2d-2}{2}
\end{array}\right\} 
\]
\[
h^{2,0}\left(V_{d_{0}}\cap V_{d}\right)=\left\{ \begin{array}{c}
dd_{0}\text{·}l_{2}^{2,0}+\left(d\left(\frac{d_{0}^{2}-d_{0}}{2}\right)+d_{0}\left(\frac{d^{2}-d}{2}\right)\right)\text{·}l_{1}^{1,0}\\
+\left(\frac{d_{0}^{2}d^{2}}{4}+\frac{dd_{0}^{3}+d_{0}d^{3}}{6}-3\frac{dd_{0}^{2}+d_{0}d^{2}}{4}+11\frac{d_{0}d}{12}\right)\text{·}l_{0}^{0,0}\\
-\frac{\left(d_{0}d^{2}+d_{0}^{2}d\right)}{2}+2d_{0}d-1
\end{array}\right\} 
\]

\[
h^{2,0}\left(V_{d_{0}}\cap V_{d_{1}}\right)=\left\{ \begin{array}{c}
d_{0}d_{1}\text{·}l_{2}^{2,0}+\left(\left(\frac{d_{0}d_{1}^{2}+d_{1}d_{0}^{2}}{2}\right)-d_{0}d_{1}\right)\text{·}l_{1}^{1,0}\\
+\left(\frac{3d_{0}^{2}d_{1}^{2}-9\left(d_{0}^{2}d_{1}+d_{0}d_{1}^{2}\right)+2\left(d_{0}d_{1}^{3}+d_{1}d_{0}^{3}\right)+11d_{0}d_{1}}{12}\right)\text{·}l_{0}^{0,0}\\
-\left(\frac{d_{1}d_{0}^{2}+d_{0}d_{1}^{2}}{2}\right)+2d_{0}d_{1}-1
\end{array}\right\} .
\]

Also for $d>1$
\[
h^{1,1}\left(V_{d_{0}}\cap V_{d}\right)=\left\{ \begin{array}{c}
h^{0,1}\left(V_{d_{0}}\cap V_{d-1}\cap V_{1}\right)+\\
h_{prim}^{1,1}\left(V_{d_{0}}\cap V_{d-1}\right)+h_{prim}^{1,1}\left(V_{d_{0}}\cap V_{1}\right)+h^{0,0}\left(V_{d_{0}}\cap V_{d-1}\cap V_{1}\cap V_{d}\right)\\
+h^{1,0}\left(V_{d_{0}}\cap V_{d-1}\cap V_{1}\right)
\end{array}\right\} 
\]
\[
h^{1,1}\left(V_{d_{0}}\cap V_{d}\right)=\left\{ \begin{array}{c}
2\left(\left(d_{0}\left(d-1\right)\right)l_{1}^{1,0}+\left(\frac{d_{0}\left(d-1\right)\left(d_{0}+d-3\right)}{2}\right)\text{·}l_{0}^{0,0}-d_{0}\left(d-1\right)+1\right)\\
+h_{prim}^{1,1}\left(V_{d_{0}}\cap V_{d-1}\right)\\
+\left(d_{0}\text{·}l_{2}^{1,1}+\left(d_{0}^{2}-d_{0}\right)l_{1}^{1,0}+\left(\frac{2d_{0}^{3}-3d_{0}^{2}+d_{0}}{3}\right)\text{·}l_{0}^{0,0}-d_{0}^{2}+d_{0}-1\right)\\
+d_{0}\left(d-1\right)d\text{·}l_{0}^{0,0}
\end{array}\right\} 
\]
\[
h^{1,1}\left(V_{d_{0}}\cap V_{d}\right)=\left\{ \begin{array}{c}
\left(2d_{0}\left(d-1\right)\right)l_{1}^{1,0}+d_{0}\left(d_{0}\left(d-1\right)+\left(d-1\right)^{2}-2\left(d-1\right)\right)\text{·}l_{0}^{0,0}\\
-2d_{0}\left(d-1\right)+2h_{prim}^{1,1}\left(V_{d_{0}}\cap V_{d-1}\right)\\
+\left(d_{0}\text{·}l_{2}^{1,1}+\left(d_{0}^{2}-d_{0}\right)l_{1}^{1,0}+\left(\frac{2d_{0}^{3}-3d_{0}^{2}+d_{0}}{3}\right)\text{·}l_{0}^{0,0}-d_{0}^{2}+d_{0}-1\right)\\
+d_{0}\left(\left(d-1\right)^{2}+\left(d-1\right)\right)\text{·}l_{0}^{0,0}
\end{array}\right\} 
\]
\[
h^{1,1}\left(V_{d_{0}}\cap V_{d}\right)=\left\{ \begin{array}{c}
\left(2d_{0}\left(d-1\right)\right)l_{1}^{1,0}+d_{0}\left(d_{0}\left(d-1\right)+\left(d-1\right)^{2}-2\left(d-1\right)\right)\text{·}l_{0}^{0,0}\\
-2d_{0}\left(d-1\right)+2+h_{prim}^{1,1}\left(V_{d_{0}}\cap V_{d-1}\right)\\
+\left(d_{0}\text{·}l_{2}^{1,1}+\left(d_{0}^{2}-d_{0}\right)l_{1}^{1,0}+\left(\frac{2d_{0}^{3}-3d_{0}^{2}+d_{0}}{3}\right)\text{·}l_{0}^{0,0}-d_{0}^{2}+d_{0}-1\right)\\
+d_{0}\left(\left(d-1\right)^{2}+\left(d-1\right)\right)\text{·}l_{0}^{0,0}
\end{array}\right\} 
\]
\[
h^{1,1}\left(V_{d_{0}}\cap V_{d}\right)=\left\{ \begin{array}{c}
\left(2d_{0}\left(\sum_{k=1}^{d-1}k\right)\right)\text{·}l_{1}^{1,0}\\
+d_{0}\left(d_{0}\left(\sum_{k=1}^{d-1}k\right)+\left(\sum_{k=1}^{d-1}k^{2}\right)-2\left(\sum_{k=1}^{d-1}k\right)\right)\text{·}l_{0}^{0,0}\\
-2d_{0}\left(\sum_{k=1}^{d-1}k\right)+2\left(d-1\right)\\
+d\left(d_{0}\text{·}l_{2}^{1,1}+\left(d_{0}^{2}-d_{0}\right)l_{1}^{1,0}+\left(\frac{2d_{0}^{3}-3d_{0}^{2}+d_{0}}{3}\right)\text{·}l_{0}^{0,0}-d_{0}^{2}+d_{0}-1\right)\\
+d_{0}\left(\left(\sum_{k=1}^{d-1}k^{2}\right)+\left(\sum_{k=1}^{d-1}k\right)\right)\text{·}l_{0}^{0,0}-\left(d-2\right)
\end{array}\right\} 
\]
\[
h^{1,1}\left(V_{d_{0}}\cap V_{d}\right)=\left\{ \begin{array}{c}
2d_{0}\left(\frac{d^{2}-d}{2}\right)l_{1}^{1,0}+\\
d_{0}\left(d_{0}\left(\frac{d^{2}-d}{2}\right)+\frac{2\left(d-1\right)^{3}+3\left(d-1\right)^{2}+\left(d-1\right)}{6}-2\left(\frac{d^{2}-d}{2}\right)\right)\text{·}l_{0}^{0,0}\\
-2d_{0}\left(\frac{d^{2}-d}{2}\right)+2\left(d-1\right)\\
+\left(d_{0}d\text{·}l_{2}^{1,1}+\left(d_{0}^{2}-d_{0}\right)d\text{·}l_{1}^{1,0}+\left(\frac{2d_{0}^{3}d-3d_{0}^{2}d+d_{0}d}{3}\right)\text{·}l_{0}^{0,0}-d_{0}^{2}d+d_{0}d-d\right)\\
-\left(d-2\right)+d_{0}\left(\frac{2d^{3}-3d^{2}+d}{6}+\frac{d^{2}-d}{2}\right)\text{·}l_{0}^{0,0}
\end{array}\right\} 
\]
\[
h^{1,1}\left(V_{d_{0}}\cap V_{d}\right)=\left\{ \begin{array}{c}
2d_{0}\left(\frac{d^{2}-d}{2}\right)l_{1}^{1,0}+\\
d_{0}\left(\left(d_{0}\left(\frac{d^{2}-d}{2}\right)+\frac{2d^{3}-3d^{2}+d}{6}-2\left(\frac{d^{2}-d}{2}\right)\right)\right)\text{·}l_{0}^{0,0}-2d_{0}\left(\frac{d^{2}-d}{2}\right)+2\left(d-1\right)\\
+\left(d_{0}d\text{·}l_{2}^{1,1}+\left(d_{0}^{2}-d_{0}\right)d\text{·}l_{1}^{1,0}+\left(\frac{2d_{0}^{3}d-3d_{0}^{2}d+d_{0}d}{3}\right)\text{·}l_{0}^{0,0}-d_{0}^{2}d+d_{0}d\right)\\
-\left(2d-2\right)+d_{0}\left(\frac{2d^{3}-2d}{6}\right)\text{·}l_{0}^{0,0}
\end{array}\right\} 
\]
\[
h^{1,1}\left(V_{d_{0}}\cap V_{d}\right)=\left\{ \begin{array}{c}
d_{0}d\text{·}l_{2}^{1,1}+\left(d_{0}\left(d^{2}-d\right)+d\left(d_{0}^{2}-d_{0}\right)\right)\text{·}l_{1}^{1,0}\\
+\left(\begin{array}{c}
\frac{d_{0}^{2}d^{2}}{2}-\frac{d_{0}^{2}d}{2}+\frac{2d_{0}d^{3}-3d_{0}d^{2}+d_{0}d}{6}\\
-d_{0}d^{2}+d_{0}d+\frac{2d_{0}^{3}d-3d_{0}^{2}d+d_{0}d}{3}\\
+\left(\frac{2d_{0}d^{3}-2d_{0}d}{6}\right)
\end{array}\right)\text{·}l_{0}^{0,0}\\
-\left(d_{0}^{2}d+d_{0}d^{2}\right)+2d_{0}d
\end{array}\right\} 
\]
\[
h^{1,1}\left(V_{d_{0}}\cap V_{d}\right)=\left\{ \begin{array}{c}
d_{0}d\text{·}l_{2}^{1,1}+\left(d_{0}\left(d^{2}-d\right)+d\left(d_{0}^{2}-d_{0}\right)\right)\text{·}l_{1}^{1,0}\\
\left(\frac{d_{0}^{2}d^{2}}{2}+2\frac{\left(d_{0}^{3}d+d_{0}d^{3}\right)}{3}-3\frac{\left(d_{0}^{2}d+d_{0}d^{2}\right)}{2}+7\frac{d_{0}d}{6}\right)\text{·}l_{0}^{0,0}\\
-\left(d_{0}^{2}d+d_{0}d^{2}\right)+2d_{0}d
\end{array}\right\} .
\]

Summarizing we now have the following table of Hodge numbers:

\smallskip{}

\noindent \begin{center}
\begin{tabular}{|c|c|}
\hline 
 & $\left(V_{d_{0}}\cap V_{d_{1}}\right)\subseteq\left|L\right|\subseteq Q$\tabularnewline
\hline 
\hline 
$h^{0,0}$ & $1$\tabularnewline
\hline 
$h^{1,0}$ & $d_{0}\text{·}d_{1}\text{·}l_{2}^{1,0}$\tabularnewline
\hline 
$h^{2,0}$ & $\begin{array}{c}
d_{0}d_{1}\text{·}l_{2}^{2,0}+\left(\left(\frac{d_{0}d_{1}^{2}+d_{1}d_{0}^{2}}{2}\right)-d_{0}d_{1}\right)\text{·}l_{1}^{1,0}\\
+\left(\frac{3d_{0}^{2}d_{1}^{2}+2\left(d_{0}d_{1}^{3}+d_{1}d_{0}^{3}\right)-9\left(d_{0}^{2}d_{1}+d_{0}d_{1}^{2}\right)+11d_{0}d_{1}}{12}\right)\text{·}l_{0}^{0,0}\\
-\left(\frac{d_{1}d_{0}^{2}+d_{0}d_{1}^{2}}{2}\right)+2d_{0}d_{1}-1
\end{array}$\tabularnewline
\hline 
$h^{1,1}$ & $+\begin{array}{c}
d_{0}d_{1}\text{·}l_{2}^{1,1}+\left(\left(d_{0}d_{1}^{2}+d_{1}d_{0}^{2}\right)-2d_{0}d_{1}\right)\text{·}l_{1}^{1,0}\\
\left(\frac{3d_{0}^{2}d_{1}^{2}+4\left(d_{0}^{3}d_{1}+d_{0}d_{1}^{3}\right)-9\left(d_{0}^{2}d_{1}+d_{0}d_{1}^{2}\right)+7d_{0}d_{1}}{6}\right)\text{·}l_{0}^{0,0}\\
-\left(d_{0}^{2}d_{1}+d_{0}d_{1}^{2}\right)+2d_{0}d_{1}
\end{array}$\tabularnewline
\hline 
\end{tabular}
\par\end{center}

\smallskip{}

\subsection{Hodge numbers of the threefolds $V_{d}\subseteq\left|L\right|\subseteq Q$}

For $d=d_{1}+d_{2}$ we again study the linear degeneration
\[
\left\{ t\text{·}F_{d}-F_{d_{1}}\text{·}F_{d_{2}}=0\right\} \subseteq\left|K_{B_{3}}\right|
\]
and take a small resolution of the nodal locus of the total space
\[
\left\{ \left|\begin{array}{cc}
t & F_{d_{2}}\\
F_{d_{1}} & F_{d}
\end{array}\right|=0\right\} \subseteq\mathbb{C}\times\left|K_{B_{3}}\right|
\]
where all four entries in the above matrix are zero. We can assume
the small resolution is absorbed in $V_{d_{2}}$ with resulting space
$\tilde{V}_{d_{2}}$. 

For $p+q=k=5$
\[
h^{p,q}\left(V_{d}\right)=h^{p-1,q-1}\left(V_{d}\cap V_{d_{1}}\cap V_{d_{2}}\right)+h^{p-1,q-1}\left(V_{d_{1}}\cap V_{d_{2}}\right)
\]
so that incorporating the above chart for lower dimensions
\[
\begin{array}{c}
h^{5,0}\left(V_{d}\right)=0\\
h^{4,1}\left(V_{d}\right)=0\\
h^{3,2}\left(V_{d}\right)=h^{3,2}\left(V_{1}\right)=l_{3}^{1,0}.
\end{array}
\]
.

For $p+q=k=4$
\[
\begin{array}{c}
h^{p,q}\left(V_{d}\right)=\\
h^{p,q}\left(\ker\left(H^{p,q}\left(V_{d_{1}}\right)\oplus H^{p,q}\left(V_{d_{2}}\right)\rightarrow H^{p,q}\left(V_{d_{1}}\cap V_{d_{2}}\right)\right)\right)\\
+h^{p-1,q-1}\left(V_{d}\cap V_{d_{1}}\cap V_{d_{2}}\right)-h^{p,q}\left(V_{d_{1}}\cap V_{d_{2}}\right)\\
=h^{p,q}\left(V_{d_{1}}\right)+h^{p,q}\left(V_{d_{2}}\right)-\left\{ \begin{array}{c}
1\,\,if\,p=q=2\\
0\,\,otherwise
\end{array}\right.
\end{array}
\]
that is, if $p+q=k=4$, 
\[
h^{p,q}\left(V_{d}\right)=\left\{ \begin{array}{c}
d\text{·}l_{3}^{p,q}-\left(d-1\right)\,\,if\,p=q=2\\
d\text{·}l_{3}^{p,q}\,\,otherwise.
\end{array}\right.
\]

Again by asymptotic mixed Hodge theory as in \cite{Clemens-0} and
arranging by weights in descending order, we have 

\begin{equation}
\begin{array}{c}
h^{3,0}\left(V_{d}\right)=\left\{ \begin{array}{c}
h^{2,0}\left(V_{d_{1}}\cap V_{d_{2}}\right)+\\
h^{3,0}\left(V_{d_{1}}\right)+h^{3,0}\left(V_{d_{2}}\right)
\end{array}\right\} =\left\{ \begin{array}{c}
h^{2,0}\left(V_{d_{1}}\cap V_{d_{2}}\right)+\\
h^{3,0}\left(V_{d_{1}}\right)+h^{3,0}\left(V_{d_{2}}\right)
\end{array}\right\} \\
h^{2,1}\left(V_{d}\right)=\left\{ \begin{array}{c}
h_{prim}^{1,1}\left(V_{d_{1}}\cap V_{d_{2}}\right)+\\
h_{prim}^{2,1}\left(V_{d_{1}}\right)+h_{prim}^{2,1}\left(V_{d_{2}}\right)\\
+h^{1,0}\left(V_{d}\cap V_{d_{1}}\cap V_{d_{2}}\right)+h^{2,0}\left(V_{d_{1}}\cap V_{d_{2}}\right)
\end{array}\right\} .
\end{array}\label{eq:Hodge}
\end{equation}
Again letting $d_{1}=1$ we use
\[
h^{2,0}\left(V_{1}\cap V_{d}\right)=d\text{·}l_{2}^{2,0}+\left(\frac{d^{2}-d}{2}\right)\text{·}l_{1}^{1,0}+\left(\frac{d^{3}-3d^{2}+2d}{6}\right)\text{·}l_{0}^{0,0}-\frac{d^{2}-3d+2}{2}
\]
to obtain
\[
h^{3,0}\left(V_{d}\right)=\left\{ \begin{array}{c}
\left(d-1\right)\text{·}l_{2}^{2,0}+\left(\frac{\left(d-1\right)^{2}-\left(d-1\right)}{2}\right)\text{·}l_{1}^{1,0}\\
+\left(\frac{\left(d-1\right)^{3}-3\left(d-1\right)^{2}+2\left(d-1\right)}{6}\right)\text{·}l_{0}^{0,0}-\frac{\left(d-1\right)^{2}-3\left(d-1\right)+2}{2}\\
+l_{3}^{3,0}\left(V_{1}\right)+h^{3,0}\left(V_{d-1}\right)
\end{array}\right\} 
\]
\[
h^{3,0}\left(V_{d}\right)=\left\{ \begin{array}{c}
d\text{·}l_{3}^{3,0}+\sum_{k=1}^{d-1}\left(d-1\right)\text{·}l_{2}^{2,0}\\
+\left(\frac{\sum_{k=1}^{d-1}\left(d-1\right)^{2}-\sum_{k=1}^{d-1}\left(d-1\right)}{2}\right)\text{·}l_{1}^{1,0}\\
+\left(\frac{\sum_{k=1}^{d-1}\left(d-1\right)^{3}-3\sum_{k=1}^{d-1}\left(d-1\right)^{2}+2\sum_{k=1}^{d-1}\left(d-1\right)}{6}\right)\text{·}l_{0}^{0,0}\\
-\frac{\sum_{k=1}^{d-1}\left(d-1\right)^{2}-3\sum_{k=1}^{d-1}\left(d-1\right)+2\left(d-1\right)}{2}
\end{array}\right\} 
\]
\[
h^{3,0}\left(V_{d}\right)=\left\{ \begin{array}{c}
d\text{·}l_{3}^{3,0}+\left(\frac{d^{2}-d}{2}\right)\text{·}l_{2}^{2,0}\\
+\left(\frac{2\left(d-1\right)^{3}+3\left(d-1\right)^{2}+\left(d-1\right)-3\left(d^{2}-d\right)}{12}\right)\text{·}l_{1}^{1,0}\\
+\left(\frac{\left(d-1\right)^{4}-2\left(d-1\right)^{3}-5\left(d-1\right)^{2}-2\left(d-1\right)+4d^{2}-4d}{24}\right)\text{·}l_{0}^{0,0}\\
-\frac{2\left(d-1\right)^{3}+3\left(d-1\right)^{2}+\left(d-1\right)-9\left(d^{2}-d\right)+12\left(d-1\right)}{12}
\end{array}\right\} 
\]
\[
h^{3,0}\left(V_{d}\right)=\left\{ \begin{array}{c}
d\text{·}l_{3}^{3,0}+\left(\frac{d^{2}-d}{2}\right)\text{·}l_{2}^{2,0}\\
+\left(\frac{2d^{3}-6d^{2}+4d}{12}\right)\text{·}l_{1}^{1,0}\\
+\left(\frac{d^{4}-6d^{3}+11d^{2}-6d}{24}\right)\text{·}l_{0}^{0,0}\\
-\frac{d^{3}-6d^{2}+11d-6}{6}
\end{array}\right\} 
\]

\noindent \begin{center}
and referring in addition to
\[
h^{1,1}\left(V_{1}\cap V_{d}\right)=d\text{·}l_{2}^{1,1}+\left(d^{2}-d\right)l_{1}^{1,0}+\left(\frac{2d^{3}-3d^{2}+d}{3}\right)\text{·}l_{0}^{0,0}-d^{2}+d
\]
we obtain
\[
h^{2,1}\left(V_{d}\right)=\left\{ \begin{array}{c}
h_{prim}^{1,1}\left(V_{1}\cap V_{d-1}\right)+\\
h_{prim}^{2,1}\left(V_{1}\right)+h_{prim}^{2,1}\left(V_{d-1}\right)\\
+h^{1,0}\left(V_{d}\cap V_{d-1}\cap V_{1}\right)+h^{2,0}\left(V_{1}\cap V_{d-1}\right)
\end{array}\right\} 
\]
\[
h_{prim}^{2,1}\left(V_{d}\right)=\left\{ \begin{array}{c}
\left(d-1\right)\text{·}l_{2}^{1,1}+\left(\left(d-1\right)^{2}-\left(d-1\right)\right)l_{1}^{1,0}\\
+\left(\frac{2\left(d-1\right)^{3}-3\left(d-1\right)^{2}+\left(d-1\right)}{3}\right)\text{·}l_{0}^{0,0}-\left(d-1\right)^{2}+\left(d-1\right)-1\\
+\left(l_{3}^{2,1}-l_{3}^{1,0}\right)+h_{prim}^{2,1}\left(V_{d-1}\right)\\
+\left(d^{2}-d\right)\text{·}l_{1}^{1,0}+\left(\frac{\left(d^{2}-d\right)\left(2d-3\right)}{2}\right)\text{·}l_{0}^{0,0}-\left(d^{2}-d\right)+1\\
+\left(d-1\right)\text{·}l_{2}^{2,0}+\left(\frac{\left(d-1\right)^{2}-\left(d-1\right)}{2}\right)\text{·}l_{1}^{1,0}\\
+\left(\frac{\left(d-1\right)^{3}-3\left(d-1\right)^{2}+2\left(d-1\right)}{6}\right)\text{·}l_{0}^{0,0}-\frac{\left(d-1\right)^{2}-3\left(d-1\right)+2}{2}
\end{array}\right\} 
\]
$h_{prim}^{2,1}\left(V_{d}\right)=\left\{ \begin{array}{c}
\left(l_{3}^{2,1}-l_{3}^{1,0}\right)+h_{prim}^{2,1}\left(V_{d-1}\right)+\left(d-1\right)\text{·}l_{2}^{1,1}+\left(d-1\right)\text{·}l_{2}^{2,0}\\
+\left(\left(d-1\right)^{2}-\left(d-1\right)\right)\text{·}l_{1}^{1,0}+\left(\left(d-1\right)^{2}+\left(d-1\right)\right)l_{1}^{1,0}\\
+\left(\frac{\left(d-1\right)^{2}-\left(d-1\right)}{2}\right)\text{·}l_{1}^{1,0}\\
+\left(\frac{4\left(d-1\right)^{3}-6\left(d-1\right)^{2}+2\left(d-1\right)}{6}\right)\text{·}l_{0}^{0,0}\\
+\left(\frac{6\left(d-1\right)^{3}+3\left(d-1\right)^{2}-3\left(d-1\right)}{6}\right)\text{·}l_{0}^{0,0}\\
+\left(\frac{\left(d-1\right)^{3}-3\left(d-1\right)^{2}+2\left(d-1\right)}{6}\right)\text{·}l_{0}^{0,0}\\
-\left(\left(d-1\right)^{2}-\left(d-1\right)\right)-\left(\left(d-1\right)^{2}+\left(d-1\right)\right)\\
-\frac{\left(d-1\right)^{2}-3\left(d-1\right)+2}{2}
\end{array}\right\} $
\[
h_{prim}^{2,1}\left(V_{d}\right)=\left\{ \begin{array}{c}
\left(l_{3}^{2,1}-l_{3}^{1,0}\right)+h_{prim}^{2,1}\left(V_{d-1}\right)+\left(d-1\right)\text{·}l_{2}^{1,1}+\left(d-1\right)\text{·}l_{2}^{2,0}\\
+\left(\frac{5\left(d-1\right)^{2}-\left(d-1\right)}{2}\right)\text{·}l_{1}^{1,0}\\
+\left(\frac{11\left(d-1\right)^{3}-6\left(d-1\right)^{2}+\left(d-1\right)}{6}\right)\text{·}l_{0}^{0,0}\\
-\left(\frac{5\left(d-1\right)^{2}-3\left(d-1\right)+2}{2}\right)
\end{array}\right\} 
\]
\[
h_{prim}^{2,1}\left(V_{d}\right)=\left\{ \begin{array}{c}
d\text{·}\left(l_{3}^{2,1}-l_{3}^{1,0}\right)+\left(\sum_{k=1}^{d-1}k\right)\text{·}l_{2}^{1,1}+\left(\sum_{k=1}^{d-1}k\right)\text{·}l_{2}^{2,0}\\
+\left(\frac{5\left(\sum_{k=1}^{d-1}k^{2}\right)-\left(\sum_{k=1}^{d-1}k\right)}{2}\right)\text{·}l_{1}^{1,0}\\
+\left(\frac{11\left(\sum_{k=1}^{d-1}k^{3}\right)-6\left(\sum_{k=1}^{d-1}k^{2}\right)+\left(\sum_{k=1}^{d-1}k\right)}{6}\right)\text{·}l_{0}^{0,0}\\
-\left(\frac{5\left(\sum_{k=1}^{d-1}k^{2}\right)-3\left(\sum_{k=1}^{d-1}k\right)+2\left(d-1\right)}{2}\right)
\end{array}\right\} 
\]
\[
h_{prim}^{2,1}\left(V_{d}\right)=\left\{ \begin{array}{c}
d\text{·}\left(l_{3}^{2,1}-l_{3}^{1,0}\right)+\left(\frac{d^{2}-d}{2}\right)\text{·}l_{2}^{1,1}+\left(\frac{d^{2}-d}{2}\right)\text{·}l_{2}^{2,0}\\
+\left(\frac{5\left(\frac{2d^{3}-3d^{2}+d}{6}\right)-\left(\frac{d^{2}-d}{2}\right)}{2}\right)\text{·}l_{1}^{1,0}\\
+\left(\frac{11\left(\frac{d^{4}-2d^{3}+d^{2}}{4}\right)-6\left(\frac{2d^{3}-3d^{2}+d}{6}\right)+\left(\frac{d^{2}-d}{2}\right)}{6}\right)\text{·}l_{0}^{0,0}\\
-\left(\frac{5\left(\frac{2d^{3}-3d^{2}+d}{6}\right)-3\left(\frac{d^{2}-d}{2}\right)+2\left(d-1\right)}{2}\right)
\end{array}\right\} 
\]
\[
h_{prim}^{2,1}\left(V_{d}\right)=\left\{ \begin{array}{c}
d\text{·}\left(l_{3}^{2,1}-l_{3}^{1,0}\right)+\left(\frac{d^{2}-d}{2}\right)\text{·}l_{2}^{1,1}+\left(\frac{d^{2}-d}{2}\right)\text{·}l_{2}^{2,0}\\
+\left(\frac{5d^{3}-9d^{2}+4d}{6}\right)\text{·}l_{1}^{1,0}\\
+\left(\frac{\left(\frac{11d^{4}-22d^{3}+11d^{2}}{4}\right)-\left(\frac{8d^{3}-14d^{2}+6d}{4}\right)}{6}\right)\text{·}l_{0}^{0,0}\\
-\left(\frac{\left(\frac{10d^{3}-15d^{2}+5d}{6}\right)-\left(\frac{9d^{2}-9d}{6}\right)+\frac{12d-12}{6}}{2}\right)
\end{array}\right\} 
\]
\[
h_{prim}^{2,1}\left(V_{d}\right)=\left\{ \begin{array}{c}
d\left(l_{3}^{2,1}-l_{3}^{1,0}\right)+\left(\frac{d^{2}-d}{2}\right)\text{·}l_{2}^{1,1}+\left(\frac{d^{2}-d}{2}\right)\text{·}l_{2}^{2,0}\\
+\left(\frac{5d^{3}-9d^{2}+4d}{6}\right)\text{·}l_{1}^{1,0}\\
+\left(\frac{11d^{4}-30d^{3}+25d^{2}-6d}{24}\right)\text{·}l_{0}^{0,0}\\
-\left(\frac{5d^{3}-12d^{2}+13d-6}{6}\right)
\end{array}\right\} .
\]
\par\end{center}

Summarizing we have the following table of Hodge numbers if $\dim B=3$:
\begin{center}
\smallskip{}
\begin{tabular}{|c|c|}
\hline 
 & $V_{d}\subseteq\left|L\right|\subseteq Q$\tabularnewline
\hline 
\hline 
$h^{0,0}$ & $1$\tabularnewline
\hline 
$h^{1,0}$ & $l_{3}^{1,0}$\tabularnewline
\hline 
$h^{2,0}$ & $d\text{·}l_{3}^{2,0}$\tabularnewline
\hline 
$h^{1,1}$ & $d\text{·}l_{3}^{1,1}-\left(d-1\right)$\tabularnewline
\hline 
$h^{3,0}$ & $\begin{array}{c}
d\text{·}l_{3}^{3,0}+\left(\frac{d^{2}-d}{2}\right)\text{·}l_{2}^{2,0}\\
+\left(\frac{d^{3}-3d^{2}+2d}{6}\right)\text{·}l_{1}^{1,0}\\
+\left(\frac{d^{4}-6d^{3}+11d^{2}-6d}{24}\right)\text{·}l_{0}^{0,0}\\
-\frac{d^{3}-6d^{2}+11d-6}{6}
\end{array}$\tabularnewline
\hline 
$h^{2,1}$ & $\begin{array}{c}
d\left(l_{3}^{2,1}-l_{3}^{1,0}\right)+l_{3}^{1,0}+\left(\frac{d^{2}-d}{2}\right)\text{·}l_{2}^{2,0}\\
+\left(\frac{d^{2}-d}{2}\right)\text{·}l_{2}^{1,1}+\left(\frac{5d^{3}-9d^{2}+4d}{6}\right)\text{·}l_{1}^{1,0}\\
+\left(\frac{11d^{4}-30d^{3}+25d^{2}-6d}{24}\right)\text{·}l_{0}^{0,0}\\
-\left(\frac{5d^{3}-12d^{2}+13d-6}{6}\right)
\end{array}$\tabularnewline
\hline 
\end{tabular}
\par\end{center}

\begin{center}
\smallskip{}
\par\end{center}

\end{document}